\documentclass[a4paper]{article}
\usepackage{epsfig}
\usepackage{amsmath}
\usepackage{amsthm}
\usepackage{amssymb}
\usepackage{enumerate}
\usepackage{ifthen}
\usepackage{graphics}
\usepackage{calc}
\usepackage{tikz}

\def\sDE{X,Y}
\newcommand{\haha}{densely dashed}
\usepackage{verbatim}

\numberwithin{equation}{section}
\newtheorem{theorem}{Theorem}[section]
\newtheorem{lemma}[theorem]{Lemma}
\newtheorem{observ}[theorem]{Observation}
\newtheorem{kobserv}[theorem]{Key Lemma}
\newtheorem{corollary}[theorem]{Corollary}
\theoremstyle{definition}
\newtheorem{definition}[theorem]{Definition}
\theoremstyle{remark}
\newtheorem{remark}[theorem]{Remark}
\def\qed{{\hfill$\Box$}}
\def\Proof{\noindent {\bf Proof}:\ }
\def\G{\mathbb{G}}
\def\S{\mathbb{S}}
\newcommand{\proba}[1]{#1}

\newif\ifdeveloping
\ifdeveloping
    \usepackage[notref,notcite]{showkeys}
\fi

\newif\ifproba
\ifproba
    \renewcommand{\proba}[1]{\ }
\fi

\newcommand\old[1]{}
\def\B{\mathbb{B}}
\newcommand{\mc}[1]{\mathcal{#1}}
\newcommand{\mbb}[1]{\mathbb{#1}}
\newcommand{\mb}[1]{\mathbf{#1}}
\newcommand{\mf}[1]{\mathfrak{#1}}

\newcommand{\inv}[1]{\overleftarrow{#1}}

\newcommand{\green}{{\mathbf{green}}}
\newcommand{\red}{{\mathbf{red}}}
\newcommand{\tope}{{\mbb T}}
\newcommand{\iter}{{\mb T }}

\newcommand{\topp}[1]{{\rm t}({#1})}
\newcommand{\bott}[1]{{\rm b}({#1})}
\newcommand{\verti}[2]{{v}_{#1}({#2})}
\newcommand{\first}[2]{\mathrm{first}_{#1}(#2)}
\newcommand{\last}[2]{\mathrm{last}_{#1}(#2)}
\newcommand{\est}{E^*}

\usetikzlibrary{backgrounds}
\newlength{\yy}
\newcounter{yyy}
\newcounter{yyyy}

\newcommand{\alapalap}[1]
    {
    \setlength{\yy}{#1  cm + 1.2 cm}
    \draw[step=1cm] (0.8,0.8) grid (\yy,\yy);
    }

\newcommand{\alap}[1]
    {
    \setlength{\yy}{#1  cm + 1.2 cm}
    \draw[step=1cm] (0.8,0.8) grid (\yy,\yy);
    \begin{pgfonlayer}{background}
    \setcounter{yyyy}{#1}
    \addtocounter{yyyy}{-1}
    \setcounter{yyy}{#1}
    \addtocounter{yyy}{+1}
        \foreach \x in {1,2,...,#1}
            {\draw (\x cm,\x cm) +(0.5,0.5)  node[fill=gray!30!white] (v) {\tiny  1};
            \draw (\x cm,0 cm) +(0.5,0.5)  node[] (v) {\tiny  \x}  ;
            \draw (0 cm,\x cm) +(0.5,0.5)  node[] (v) {\tiny  \x}  ;
            };
        \foreach \x in {1,2,...,\value{yyyy} }
            {\draw (\x cm,\x cm) +(0.5,1.5)  node[fill=gray!30!white] (v) {\tiny  0};
            };
    \rhh{#1}{1}{0}
    \end{pgfonlayer}
    }

\newcommand{\kalap}[1]
    {
    \setlength{\yy}{#1  cm + 1.2 cm}
    \draw[step=1cm] (0.8,0.8) grid (\yy,\yy);
    \begin{pgfonlayer}{background}
    \setcounter{yyyy}{#1}
    \addtocounter{yyyy}{-1}
    \setcounter{yyy}{#1}
    \addtocounter{yyy}{+1}
        \foreach \x in {1,2,...,#1}
            {\draw (\x cm,\x cm) +(0.5,0.5)  node[fill=gray!30!white] (v) {\tiny 1};
             };
        \foreach \x in {1,2,...,\value{yyyy} }
            {\draw (\x cm,\x cm) +(0.5,1.5)  node[fill=gray!30!white] (v) {\tiny  0};
            };
    \end{pgfonlayer}
    }

\newcommand{\hh}[3]{
    \draw (#1 cm,#2 cm) +(0.5,0.5)  node[fill=white]  {\tiny #3}  ;
    }

\newcommand{\rhh}[3]{
    \draw (#1 cm,#2 cm) +(0.5,0.5)  node[fill=gray!30!white]  {\tiny #3}  ;
    }

\newcommand{\hhc}[3]{
    \draw (#1 cm,#2 cm) +(0.5,0.5)  node[fill=white]  {\tiny #3}  ;
    \draw (#1 cm,#2 cm) +(0.5,0.5) circle (0.4cm);
    }

\newcommand{\rhhc}[3]{
    \draw (#1 cm,#2 cm) +(0.5,0.5)  node[fill=gray!30!white]  { \tiny #3}  ;
    \draw (#1 cm,#2 cm) +(0.5,0.5)	circle (0.4cm);
    }

\def\<{\left\langle}
\def\>{\right\rangle}

\title{Towards  random uniform sampling of bipartite graphs with given degree
sequence\footnotemark[1] \footnotemark[2]\footnotetext[1]{This research was  supported in part by the Hungarian Bioinformatics MTKD-CT-2006-042794, Marie Curie Host Fellowships for Transfer of Knowledge.}
\footnotetext[2]{IM and PLE  acknowledge financial support from grant \#FA9550-12-1-0405 from the U.S. Air Force Office of Scientific Research (AFOSR) and the Defense Advanced Research Projects Agency (DARPA).}}
\author{
    Istv\'an Mikl\'os\thanks{Partly supported by Hungarian NSF, under contract Nos.\ NK 78439 and PD84297} \quad  and\quad
    P{\'e}ter L. Erd{\H o}s\thanks{Partly supported by Hungarian NSF, under contract Nos.\ NK 78439 and K68262}
    \quad  and\quad
Lajos Soukup\thanks{Partly supported by Hungarian NSF, under contract Nos.\ NK  83726 and K68262}\\
    Alfr\'ed R{\'e}nyi Institute of Mathematics,\\
    Hungarian Academy of Sciences,\\
    Budapest, P.O.~Box 127, H-1364 Hungary\\
    {\texttt{ <miklosi,elp,soukup>@renyi.hu}}
        }
\begin{document}
\maketitle
\begin{abstract}
In this paper we consider a simple Markov chain for bipartite graphs with given
degree sequence on $n$ vertices.
We show that the mixing time of this Markov chain is bounded above by
a polynomial in $n$ in case of {\em half-regular} degree sequence. The
novelty of our approach lies in the construction of the multicommodity
flow in Sinclair's method.
\end{abstract}

\section{Introduction}\label{sec:intro}
\noindent
The {\em degree sequence}, $d(G)$, of a graph $G$ is the  non-increasing
sequence of its vertex degrees. A  sequence $\mathbf{d}= (d_1,\ldots,d_n)$ is
{\em graphical} iff $d(G) = \mathbf{d}$ for some simple graph $G$, and $G$ is a
{\em graphical realization} of $\mathbf{d}$.

Already at the beginning of the systematic graph theoretical research (late
fifties and early sixties) there were serious efforts to decide whether a
non-increasing sequence is graphical. Erd\H{o}s and Gallai (1960, \cite{EG60})
gave a necessary and sufficient condition, while Havel (1955, \cite{H55}) and
Hakimi (1962, \cite{H62}) independently developed a greedy algorithm to built  a
graphical realization  if there exists any. (For  more details   see for example
\cite{KTEMS}.)

Generating some (or all possible) graphs realizing a given degree sequence or
finding a typical one among the different realizations are ubiquitous problems
in network modeling, ranging from social sciences to chemical compounds and
biochemical reaction networks in the cell. (See for example the book
\cite{NBW06} for a detailed analysis, or the paper \cite{KTEMS} for a short
explanation.)

When the  number of different realizations is small, then the uniform sampling
of the different realizations can be carried out by generating all possible ones
and choosing among them uniformly.

However in cases where there are many different realizations this approach can not work. In these cases  some stochastic processes can provide solutions. Here we mention only one  of the preceding results: Molloy and Reed (1995, \cite{MR95}) applied the {\em configuration   model}  (Bollob\'as (1980, \cite{Bo80}) for the problem. (In fact, Wormald had used   it already in 1984  to generate random regular graphs of {\em moderate} degrees \cite{Wo84}.) They successfully used the model to generate random graphs with given degree sequences where the degrees are (universally)  bounded. It is well known that this method is computationally infeasible in case of general, unbounded degree sequences.

A different method was proposed by Kannan, Tetali and Vempala (1995,
\cite{KTV97}), which is based on the powerful Metropolis-Hastings algorithm:
some local transformation generates a random walk on the  family  of all
realizations. They conjectured that this process is {\em rapidly mixing} i.e.
starting from an arbitrary realization of the degree sequence the process
reaches a completely random realization in reasonable (i.e. polynomial) time.
However,  they could prove it only for bipartite regular graphs. Their
conjecture was proved for arbitrary regular graphs by Cooper, Dyer and Greenhill
(2007, \cite{CDG07}).

The original  goal of this paper  was to attack Kannan, Tetali and Vempala's conjecture for arbitrary bipartite degree sequences, performing a more subtle choice of {\em multicommodity flow}. We obtained the following result:
\begin{theorem}\label{th:main}
The Markov process - defined by Kannan, Tetali and Vempala - is  rapidly mixing  on
each bipartite half-regular degree sequence. $($In these bipartite graphs the  degrees in one vertex class are constant.$)$
\end{theorem}
\noindent Actually, we achieved somewhat more: our construction method can be used as a plug-in to a more advanced method for general degree sequences: if two particular graphical realizations at hand differ in edges which can be partitioned into alternating cycles, such that no cycle contains a chord which is an edge of another cycle in the partition, then our {\em friendly path} method provides a good multicommodity flow.

\section{Basic definitions and preliminaries}\label{sec:def}
Let $G=\left(U, V;E\right )$ be a simple bipartite graph (no parallel edges) with vertex classes $U=\{ u_1,\ldots,u_k\}, V=\{v_1,\ldots,v_l\}$. The {\em (bipartite) degree sequence} of $G$, $\mathrm{bd}(G)$ is defined as follows:
$$
\mathrm{bd}(G)=\Big(\big (d(u_1), \ldots,d(u_k)\bigr ), \bigl (d(v_1),\ldots ,d(v_l)\bigr
)\Big ),
$$
where the vertices are ordered such that both sequences are non-increasing. From now on when we say ``degree sequence'' of a bipartite graph, we will always mean the bipartite degree sequence. We will use $n$ to denote the number of vertices, that is $n=k+l.$

A pair $ (\mathbf{a}, \mathbf{b})$ of sequences is a {\em $($bipartite$)$ graphical sequence} (BGS for short) if $ (\mathbf{a}, \mathbf{b})=\mathrm{bd}\left (G \right )$ for some simple bipartite graph $G$, while the graph $G$ is a {\em (graphical) realization} of $ (\mathbf{a}, \mathbf{b})$.

Next we define the swaps, our basic operation on bipartite graphs.
\begin{definition}
Let $G=(U,V;E)$ be a bipartite graph, $u_1, u_2 \in U$, $v_1, v_2,\in V$, such
that induced subgraph $G[u_1,u_2; v_1, v_2]$ is a 1-factor, (i.e. $(u_1, v_j), (u_2,
v_{3-j})$ $ \in E$, but $(u_1, v_{3-j}), (u_2, v_j) \notin E$ for some $j$.) Then we
say that the {\bf swap on  $(u_1, u_2; v_1, v_2)$} is {\bf allowed}, and it
transforms the graph $G$ into a graph $G'=(U,V;E')$ by replacing the edges
$(u_1, v_j), (u_2, v_{3-j})$  by edges  $(u_1, v_{3-j})$ and  $(u_2, v_j)$, i.e.
\begin{equation}
E'=E\setminus \{(u_1, v_j), (u_2, v_{3-j})\} \cup \{(u_1, v_{3-j}), (u_2,
v_j)\}.
\end{equation}
\end{definition}
\noindent So a swap transforms one realization of the BGS to another (bipartite
graph) realization of the same BGS. The following proposition is a classical
result of Ryser (1957, \cite{R57}).

\begin{theorem}[Ryser]\label{th:ryser}
Let  $G_1=\left(U, V;E_1\right )$  and $G_2=\left(U, V;E_2\right )$ be two
realizations of the same BGS. Then there exists a sequence of swaps which
transforms $G_1$ into $G_2$ through different realizations of the same BGS.
\end{theorem}
\noindent Ryser's result used the language of 0 - 1 matrices. Here, to make the paper self contained, we give a short proof, using the notion of swaps. The proof is based on a well known observation of Havel and Hakimi (\cite{H55,H62}):

\begin{lemma}[Havel and Hakimi]\label{th:HH} Let $G=\left(U, V;E\right )$ be a simple bipartite graph, and assume that $d(u') \le d(u)$, furthermore $(u',v) \in E$ and $ (u,v)\not\in E.$ Then there exists a vertex $v'$ such that the swap on ${(u,u';v,v')}$ is allowed, and so  it produces a bipartite graph $G'$ from $G$ such that $\Gamma_{G'}(v) = (\Gamma_G(v) \setminus \{ u'\}) \cup \{ u \},$ where, as usual, $\Gamma_G(v)$ is the set of neighbors of $v$ in $G.$
\end{lemma}
\Proof By the pigeonhole principle there exists a vertex $v' \ne v$ such that
$(u,v') \in E$ and $ (u',v') \not\in E.$ So the swap defined on vertices $(u,u';v,v')$ is allowed. \qed

\medskip\noindent We say that the previous operation is {\em pushing up} the
neighbors of vertex $v$. Applying the pushing up operation $d$ times we obtain
the following push up lemma.
\begin{lemma}[Havel and Hakimi]\label{lm:push}
If   $G=\left(U, V;E\right )$ is a simple bipartite graph, $d(u_1)\ge d(u_2)\ge \dots\ge d(u_k)$ and $v\in V$, $d=d(v)$.  Then there is a sequence $S$  of  $d$ many swaps which transforms $G$  into a graph $G'$ such that $\Gamma_{G'}(v)=\{u_1,\dots, u_d\}$.
\end{lemma}
\noindent This pushing-up lemma also suggests (and proves the correctness of) a
greedy algorithm to construct a concrete realization of a BGS  $\big
(\mathbf{a}, \mathbf{b}\big )$.

\medskip\noindent{\bf Proof of Theorem~\ref{th:ryser}}: We prove the following
stronger statement:
\begin{itemize}\label{obs:e}
\item[(\maltese)] {\em there exists a sequence of $2e$  swaps which transforms
$G_1$ into $G_2$, where $e$ is the number of edges of $G_i$.}
\end{itemize}
We will show that any particular realization can be transformed into the same {\em canonical realization} with at most $e$ swaps. We will do it recursively: taking one by one the vertices $v_1,v_2,\ldots,v_l$ from $V$ we will define their neighbors in $U.$ After every step of the process we update the remaining degree sequence of $U$, and reorder its actual content.

To do so we introduce the following lexicographic order on the actual remaining $d(u)$ degree sequence. We always take them non-increasing order, and whenever two vertices have the same actual degree, then we take first the vertex with bigger subscript.

So take $v_1$ and by multiple applications of the Push-up Lemma \ref{lm:push} there is a sequence $T_1$ of at most $d=d(v_1)$ many swaps which transforms
$G_1$ into  a $G'_1$  such that $\Gamma_{G_1'}(v_1)=\{u_1,\dots, u_d\}$ (The actually required push up operations can be smaller if some of the first $d$ vertices were originally adjacent to $v_1$.)

We consider the bipartite graphs $G_1''=G_1'\setminus \{v_1\}$ i.e. we remove the vertex $v_1$ and all the edges connected to $v_1$. Now we reorder the vertices in the actual $U$ according to our lexicographic order, and repeat the recursive operation.

In this way after at most $\sum_{i=1}^l d(v_i)=e$ swaps we transformed $G_1$ into a well defined canonical realization $R$, furthermore this $R$ is independent from the original realization.

Now we can easily finish the proof of Theorem \ref{th:ryser} observing that if a swap transforms $H$ into $H'$, then the ``inverse swap" (choosing the same four vertices, and changing back the edges) transforms $H'$ into $H$. So if the swap sequence $T_1$ transforms $G_1$ into $R$ then it has an inverse swap-sequence $T_1'$ which transforms $R$ into $G_1$. \qed

\medskip\noindent We use this upper bound for convenience: for us a linear upper bound on this value is enough to show the polynomial upper bound of the sampling process. If somebody wanted to get tight (or at least better) upper bounds on the sampling process, then a better estimation is necessary for the swap-distance. Recently it was shown that the swap-distance $\mathbf{dist}(G_1,G_2)$ for any two realizations is smaller than
$$
\Delta := \frac{1}{2} \big |E(G_1) \Delta E(G_2) \big |.
$$
In the forthcoming paper \cite{EKM} a formula for $\mathbf{dist}(G_1,G_2)$  is determined: this is in the form of $\Delta - \alpha$ where the parameter $\alpha \ge 1.$ Unfortunately the parameter $\alpha$ is hard to determine.

\section{The Markov chain $(\G,P)$}
For a bipartite graphical sequence $\big (\mathbf{a}, \mathbf{b}\big )$ (on the fixed vertex bipartition $(U,V)$) - following Kannan, Tetali and Vempala's lead - we define a Markov chain $(\G,P)$ in the following way. $\G$ is a graph, the vertex set $V( \G)$ of the graph $\G$ consists of all possible realizations of our BGS, while the edges represent the possible swap operations: two realizations are connected if there is a swap operation which transforms one realization into the other one (and, recall, the inverse swap transforms the second one to the first one as well).

\noindent Let $P$ denote the {\em transition matrix}, which is defined as follows: if the current realization (state of the process) is $G$ then with probability $\frac{1}{2}$ we stay in the
current state (namely, we define a lazy Markov chain) and with probability
$\frac{1}{2}$ we choose uniformly two-two vertices $u_1,u_2;v_1,v_2$ from
classes $U$ and $V$ respectively and perform the swap if it is possible and move
to $G'$. Otherwise we do not perform a move. The swap moving from $G$ to $G'$ is
unique, therefore the probability of this transformation (the {\em jumping
probability} from $G$ to $G'\ne G$) is:
\begin{equation}\label{eq:prob}
\mathrm{Prob}(G\rightarrow G'):= P(G' | G) =
\frac{1}{2\genfrac{(}{)}{0pt}{}{k}{2} \genfrac{(}{)}{0pt}{}{l}{2}}.
\end{equation}
The probability of transforming $G'$ to $G$ is time-independent.
 The transition probabilities are {\em time} and edge {\em
  independent} and they are also {\em symmetric}. Therefore
$P$ is a symmetric matrix, where all off-diagonal, non-zero elements are the same, while the entries in
the main-diagonal are non-zero, but (probably) different values.

We use the convention that upper case letters $X,Y$ and $Z$ stands for vertices
of $V(\G)$.

The graph $\G$ clearly may have exponentially many vertices (that many different
realizations of the degree sequence). However, by the statement (\maltese)  (in
the proof of Theorem~\ref{th:ryser}), its diameter  is always relatively small:

\begin{corollary}\label{lm:kicsi}
The {\em swap distance} of any two realizations is at most $2e$, where $e$ is
the number of edges.
\end{corollary}

As we observed, the  graph $\G$ is connected, therefore the Markov process is
{\em irreducible}. Since our Markov chain is lazy, it is clearly aperiodic.
Finally since, as we saw, the jumping probabilities are symmetric, that is
$P(G|G') = P(G'|G)$, therefore our lazy Markov process is reversible with the
uniform distribution as the globally stable stationary distribution.

\section{Sinclair's Method}\label{sec:sinclair}
To start with we recall some definitions and notations from the literature. Since our Markov chain converges to the uniform distribution, we write all theorems for the special uniform distribution case even if the theorem holds for more general distribution, to simplify the notations. Let $P^t$ denote the $t$th power of the transition probability matrix and define
\begin{displaymath}
\Delta_X(t):=\frac 12\sum_{Y\in V(\G)} \left | P^t(Y | X)- 1/N \right |,
\end{displaymath}
where $X$ is an element of the state space of the Markov chain and $N$
is the size of the state space. We define the {\em mixing   time} as
\begin{displaymath}
{\tau}_X({\varepsilon}):=\min_{t}\big \{\Delta_X(t')\le {\varepsilon}
\text{ for all $t'\ge t$} \big \}.
\end{displaymath}
Our Markov chain is said to be {\em rapidly mixing} iff
\begin{displaymath}
{\tau}_X({\varepsilon})\le O\Big (\mathrm{poly}\big (log( N/{\varepsilon}) \big
) \Big )
\end{displaymath}
for any $X$ in the state space. Consider the different eigenvalues of $P$ in non-increasing order:
\begin{displaymath}
1={\lambda}_1>{\lambda}_2\ge \dots \ge{\lambda}_N \ge -1.
\end{displaymath}
The {\em relaxation time} ${\tau}_{rel}$ is defined as
\begin{displaymath}
{\tau}_{rel}=\frac 1 {1-{\lambda}^*}
\end{displaymath}
where $\lambda^*$ is the {\em second largest eigenvalue modulus},
\begin{displaymath}
\lambda^*:= \max\{ \lambda_2 ,|\lambda_N|\}.
\end{displaymath}
However, the eigenvalues of any lazy Markov chain are non-negative, so we do know that $\lambda^* = \lambda_2$ for our Markov chain. The following result was proved implicitly by Diaconis and Strook in 1991, and explicitly stated by Sinclair: \cite[Theorem 5']{S92}
\begin{theorem}[Sinclair] \label{tm:DS}
\qquad $\displaystyle {\tau}_x({\varepsilon})\le {\tau}_{rel}\cdot \mathrm{poly}
\big (log( N/{\varepsilon}) \big ).
$ \qed
\end{theorem}
\noindent So one way to prove that our Markov chain is
rapidly mixing
is to find a polynomial  upper
bound on ${\tau}_{rel}$. We need rapid convergence of the process to the
stationary distribution otherwise the method cannot be used in practice.

Kannan, Tetali and Vempala in \cite{KTV97} could prove that the relaxation time
of the Markov chain $(\G,P)$ is a polynomial function of the size $n := 2k$ of
$\big (\mathbf{a}, \mathbf{b}\big )$ if it is a regular bipartite degree
sequence. Here we extend their proof to show that the process is rapidly mixing
for the half-regular bipartite case.

\bigskip\noindent There are several different methods to prove fast
convergence, here we use - similarly to \cite{KTV97} - Sinclair's {\em
  multicommodity flow method} (\cite{S92}).

\begin{theorem}\label{th:sinc}
Let $\mathbb{H}$ be a graph whose vertices represent the possible states of a
time reversible finite state Markov chain $\mathcal{M}$, and where $(U,V) \in
E(\mathbb{H})$ iff the transition probabilities of $\mathcal{M}$ satisfy $P(U|V)
P(V|U) \ne 0.$ For all $X\ne Y \in V(\mathbb{H})$ let $\Gamma_{X,Y}$ be a set of
paths in $\mathbb{H}$ connecting $X$ and $Y$ and let $\pi_{\raisebox{-1pt}{\scriptsize X,Y}}$ be a probability distribution  on $\Gamma_{X,Y}$.
Furthermore let
\begin{displaymath}
\Gamma := \bigcup_{ X\ne Y \in V(\mathbb{H})} \Gamma_{X,Y}
\end{displaymath}
where the elements of $\Gamma$ are called {\em paths}. We also assume
that there is a stationary distribution $\pi$ on the vertices
$V(\mathbb{H})$.  We define the capacity of an edge $e=(W,Z)$ as
\begin{displaymath}
Q(e) := \pi(W) P(Z|W)
\end{displaymath}
and we denote the length of a path $\gamma$ by $|\gamma|$.
Finally let
\begin{equation}\label{eq:multiflow}
\kappa_{\Gamma} := \max_{e \in E(\mathbb{H})}
\frac{1}{Q(e)}
    \sum_{\genfrac{}{}{0pt}{}{\scriptstyle  X,Y \in
    V(\mathbb{H})}{\scriptstyle  \gamma \in \Gamma_{X,Y}\ : \ e \in
    \gamma}} \pi(X) \pi(Y)\pi_{\raisebox{-1pt}{\scriptsize X,Y}}(\gamma) |\gamma|.
\end{equation}
Then
\begin{equation}\label{eq:relax1}
 \tau_{\mathrm{rel}}(\mathcal{M})  \le \kappa_\Gamma
\end{equation}
holds. \qed
\end{theorem}
\bigskip\noindent We are going to apply Theorem~\ref{th:sinc} for our Markov chain $(\G, P$).  Using the notation $|V(\G)|:=N$, the (uniform) stationary distribution has the value {$\pi(X)=1/N$} for each vertex $X\in V(\G).$ Furthermore each transition probability has the property $P(X|Y) \ge 1/ n^4$ (recall that $n=k+l$, that is $n$ denotes the number of the vertices of any realization).  So if we can design a multicommodity flow such that  each path is shorter then an appropriate $\mathrm{poly}(n)$ function, then simplifying inequality~(\ref{eq:multiflow}) we can turn inequality (\ref{eq:relax1}) to the form:
\begin{equation}\label{eq:multiflow_simple}
\tau_{\mathrm{rel}} \leq \frac{\mathrm{poly}(n)}{N} \left (\max_{e \in E(\mathbb{H})}
\sum_{\genfrac{}{}{0pt}{}{\scriptstyle X,Y \in     V(\mathbb{H})}{\scriptstyle \gamma \in \Gamma_{X,Y}\ : \ e \in     \gamma}}\pi_{\raisebox{-1pt}{\scriptsize X,Y}}(\gamma)\right ).
\end{equation}
If $Z\in e$, then
\begin{equation}\label{eq:multiflow_simplex1}
\sum_{\genfrac{}{}{0pt}{}{\scriptstyle X,Y \in     V(\mathbb{H})}{\scriptstyle
\gamma \in \Gamma_{X,Y}\ : \ e \in     \gamma}}\pi_{\raisebox{-1pt}{\scriptsize X,Y}}(\gamma) \leq
\sum_{\genfrac{}{}{0pt}{}{\scriptstyle X,Y \in     V(\mathbb{H})}{\scriptstyle
\gamma \in \Gamma_{X,Y}\ : \ Z \in     \gamma}}\pi_{\raisebox{-1pt}{\scriptsize X,Y}}(\gamma),
\end{equation}
so we have
\begin{equation}\label{eq:multiflow_simple1.5}
\tau_{\mathrm{rel}} \leq \frac{\mathrm{poly}(n)}{N} \left (\max_{Z \in
V(\mathbb{H})}
\sum_{\genfrac{}{}{0pt}{}{\scriptstyle X,Y \in     V(\mathbb{H})}{\scriptstyle
\gamma \in \Gamma_{X,Y}\ : \ Z \in     \gamma}}\pi_{\raisebox{-1pt}{\scriptsize X,Y}}(\gamma)\right ).
\end{equation}

We make one more assumption. Namely, that for each $X,Y\in V(\mbb G)$ there is a non-empty finite set $S_{X,Y}$ (which draws its elements from a pool of symbols) and for each $s\in S_{X,Y}$ there is a path $\varUpsilon(X,Y,s)$ from $X$ to $Y$ such that
\begin{equation}
\Gamma_{X,Y}= \{\varUpsilon(X,Y,s) : s\in S_{X,Y}\}.
\end{equation}
It can happen that $\varUpsilon(X,Y,s) = \varUpsilon(X,Y,s')$ for $s\ne s'$,
so we consider $\Gamma_{X,Y}$ as a ``multiset'' and so we should take
$$
\pi_{X,Y}(\gamma)=\frac{\left | \big \{s\in S_{X,Y}: \gamma=\varUpsilon(X,Y,s) \big \} \right |
}{\left |S_{X,Y}\right |}
$$
for ${\gamma}\in \Gamma_{X,Y}$.

\bigskip\noindent
Putting together the observations and simplifications above we obtain the\\
{\bf Simplified Sinclair's method:}\\
For each $X,Y\in V(\mbb G)$ find a non-empty finite set $S_{X,Y}$ and for each $s\in S_{X,Y}$ find a path $\varUpsilon(X,Y,s)$ from $X$ to $Y$ such that
\begin{itemize}
\item each path is shorter than an appropriate $\mathrm{poly}(n)$ function,
\item  for each $Z\in V(\mbb G)$
\begin{equation}\label{eq:multiflow_simple2}
\sum_{\genfrac{}{}{0pt}{}{\scriptstyle X,Y \in     V(\mathbb{G})}{ }}
\frac{\left | \big \{s\in S_{X,Y}: Z\in \varUpsilon(X,Y,s)\big \} \right | }{|S_{X,Y}|}
\\\le poly(n)\cdot N.
\end{equation}
Then our Markov  chain  $(\mbb G, P)$ is rapidly mixing.
\end{itemize}

\section{Multicommodity flow - general considerations }\label{sec:general_new}
Our construction method for multicommodity flow commences on the trail of
Kannan, Tetali and Vempala (\cite{KTV97}), and Cooper, Dyer and Greenhill
(\cite{CDG07}). However the main difference among these papers lays in the method of the construction of the multicommodity flow.

We fix a bipartite graphical sequence  $\big (\mathbf{a}, \mathbf{b}\big )$,
and consider the graph $\G$ where the vertices of  $\G$ are the realizations of $\big (\mathbf{a}, \mathbf{b}\big )$, while the edges correspond to the possible swap operations. Therefore if $X\in \G$, then $X$ is a simple bipartite graph $(U,V;E(X))$, where $U$ and $V$ are fixed finite sets.

\medskip\noindent We can outline the construction of the path system from $X\in \G$ to $Y\in \G$ as follows:
\begin{enumerate}[(Step 1)]
\item We decompose the symmetric difference $\Delta$ of   $E(X)$ and $E(Y)$
into alternating circuits:
    \begin{displaymath}
            W_1,W_2\dots, W_{k_s}.
    \end{displaymath}
The construction uses the method of \cite{CDG07} to parameterize all the possible decompositions (see Subsection \ref{ss:circle}). Roughly speaking, the parameter set $S_{X,Y}$ will be the collection of all pairing of edges $E(X) \setminus E(Y)$ and $E(Y)\setminus E(X)$ adjacent to $w$, for all $w\in U\cup V$.
\item We decompose every alternating circuit $W_i$  into alternating cycles
    \begin{displaymath}
            C^i_1,C^i_2\dots, C^i_{k_i},
    \end{displaymath}
and we will construct the canonical path from $X$ to $Y$ in such a way that first we switch the edges  $E(X)\setminus E(Y)$ and $E(Y)\setminus E(X)$ in
$C^i_1$, then  in $C^i_2$, etc.

Let's denote $Z$ an arbitrary vertex along the canonical path. To apply Sinclair method we will need that  the elements of  $S_{X,Y}$ can be reconstructed  from elements of $S_{\Delta\cap E(Z),\Delta\setminus E(Z)}$ (using another small parameter set). In \cite{CDG07} the authors could prove that the  elements of $S_{X,Y}$ are ``almost'' in $S_{\Delta\cap E(Z),\Delta\setminus E(Z)}$. Unfortunately, it is not true for our construction. This is the reason that we should introduce  a much more complicated ``reconstruction'' method
in Subsection \ref{ss:rec} below.
\end{enumerate}

\subsection{Alternating circuit decompositions}\label{ss:circle}
Before we start  this subsection we should recall some definitions:
\begin{definition}
In a simple graph, a sequence of pairwise disjoint edges $e_1,\ldots,$ $ e_t $
forms a {\bf circuit} iff there are vertices $v_1,\dots, v_t $ such that
$e_i=(v_i,v_{i+1})$  (the summation is performed modulo $t$). This circuit is a
{\bf cycle} iff the vertices $v_1,\dots, v_t $ are pairwise distinct.
\end{definition}

\medskip\noindent Now let $K=(W,F\cup F')$  be a simple graph where $F\cap
F'=\emptyset$ and assume that for each vertex $w \in W$ the $F$-degree and
$F'$-degree of $w$ are the same: $d(w)=d'(w)$ for all $w\in W$. An {\em
alternating circuit decomposition} of $(F, F')$ is a circuit decomposition
such that no two consecutive edges of any circuit are in $F$ or in $F'$. Next we
are going to parameterize the alternating circuit decompositions.

The set of all edges in $F$ (in $F'$)     which are incident to a vertex $w$ is
denoted by $F(w)$  (by $F'(w)$, respectively).

If $A$ and $B$ are sets, denote $[A,B]$ the complete bipartite graph with
classes $A$ and $B$. Let
\begin{multline}
\mbb S(F,F')=\big \{s:\text{$s$ is a function, $dom(s)=W$, and for all $w\in
W$}\\
\text{$s(w)$ is a 1-factor of the complete bipartite graph $\left [ F(w),
F'(w)\right ]$}\big \}.
\end{multline}

\medskip
\begin{lemma}\label{lm:decomp}
There is a natural one-to-one correspondence between the elements of $\S(F,F' )$
and the family of all alternating circuit decompositions of $(F,F')$.
\end{lemma}
\begin{proof}
If $\mc C=\{C_1,C_2,\dots,C_n\}$ is an alternating circuit decomposition of
$(F,F')$, then define $s_{\mc C}\in \mbb S(F,F')$ as follows:
\begin{multline}\label{eq:matching}
s_{\mc C}(w):=\big \{\big((w,u),(w,u')\big)\in [F(w), F'(w)]:\\
 \text{$(w,u)$ and $(w,u')$ are consecutive edges in some  $C_i\in \mc C$} \big
\}.
\end{multline}

\noindent
On the other hand, to each  $s\in \mbb S(F,F')$ assign  an alternating circuit
decomposition
\begin{displaymath}
   \mc C_s=  \{W^s_1,W^s_2\dots, W^s_{k_s}\}
\end{displaymath}
of $(F,F')$ as follows: Consider the bipartite graph $\mathcal{F }=\big (F,F',
R(s)\big )$, where
\begin{equation}\notag
R(s)=  \big\{ \big(  (u,w), (u',w)\big ): w\in W \text{ and }
\big(  (u,w), (u',w)\big )\in s(w) \big\}.
\end{equation}
  $\mathcal{F}$ is a $2$-regular graph because for each edge $(u,v)\in F\cup F'
$
there is exactly one  $(u,w)\in F\cup F'$ with $\big(  (u,w), (u,w)\big )\in
s(u)$,
there is exactly one  $(t,v)\in F\cup F'$ with $\big(  (u,v), (t,v)\big )\in
s(v)$,
therefore the $\mathcal{F}$-neighbors of $(u,v)$ are $(u,w)$ and $(t,v)$.

$\mathcal{F}$ is a $2$-regular, so it is the union of vertex disjoint cycles
$\{W^s_i:i\in I\}$. Now $W^s_i$ can also be viewed as a sequence of edges in
$F\cup F'$, which is  an alternating circuit in $\<W,F\cup F'\>$, so
$\{W^s_i:i\in I\}$   is an alternating circuit decomposition of $(F,F')$. Since
\begin{displaymath}
 s_{{\mc C}_{\normalsize s}}=s,
\end{displaymath}
we proved the Lemma.
\end{proof}
\medskip\noindent
If the $F$-degree sequence (and therefore the $F'$-degree sequence) is
$d_1,\dots d_k$, then write
\begin{equation}
t_{F,F'}= \prod_{i=1}^k (d_i!).
\end{equation}
Clearly
\begin{displaymath}
\left |\S \big (F,F' \big ) \right |=t_{F,F'}.
\end{displaymath}

\subsection{Cycle decompositions and circuit reconstructions}
\label{ss:rec}

\noindent In this subsection we make preparations for constructing our multicommodity flow: we describe how we decompose an alternating circuit  into alternating cycles.

The problem of this venture is the following: we know along the process the symmetric difference of the edge sets of realizations $X$ and $Y$ but we do not know the distribution of the edges among $E(X)$ and $E(Y)$. If the alternating circuit under investigation is large then its cycle decomposition can contain a linear number of alternating cycles. Each cycle consists of an even number of edges, equally distributed between $X$ and $Y$. Along the process each cycle needs a parameter representing whether that particular cycle was already processed or not (which, in turn, tells which edges belong to $X$ and $Y$). Therefore along all decompositions the set of all possible parameter values can be exponentially big, which is not suitable to prove fast mixing property. Therefore we need to find another way to deal with the reconstruction problem. We can proceed as follows:

\medskip\noindent
Let $\mb x=(x_1,x_2,\dots, x_m)$ be a sequence, then we write $\inv{\mb x} =(x_m,\dots, x_2,x_1)$ for the oppositely ordered  sequence. (Here we consider $\inv{. }$ to be an operator.)

Assume that $K=(W,F\cup F')$  is  a simple bipartite graph where $F\cap
F'=\emptyset$ (in our applications we have $|F| = |F'|$), and the sequence
$$
 \mb e=(e_1,e_2,\dots,e_m)
$$
of edges is an alternating  walk in $K$. (Consequently the elements of the sequence are pairwise distinct.) In this subsection we also use extensively the notation $\mb e= e_1e_2\cdots e_m$ for the same sequence. When $\mb f $ and $\mb g$  are two sequences, then $\mb f \mb g$ stands for their concatenation. We will write $e_i=v_iv_{i+1}$. We will also use $\bott{e_i}=v_i$ (for the {\em bottom} of the edge) and $\topp{e_i}=v_{i+1}$ (for the {\em top} of the edge, considering the actual orientations along the walk). So  $\bott{e_{i+1}}=\topp{e_i}$, and $\topp{e_m}=\bott{e_1}$ iff $\mb e$ is a  circuit.

We will use the notations $\mb e(i)=e_i$ and $\verti{\mb e}j=v_j$ for $1\le i\le m$ and $1\le j\le m+1$ (the $i$th edge and the $j$th vertex of the walk). If $\mb c=e_i\cdots e_j$ is a consecutive subsequence of $\mb e$, we will also write $\bott{\mb c}=\bott{e_i}$ and $\topp{\mb c}=\topp{e_j}$. Finally $\first{}{\mb f}$ denotes the first edge,  while $\last{}{\mb f}$ denotes the last edge of walk $\mb f.$

Now let $f$ be a coloration of the edges along the walk  $f:\mb e\to \{\green, \red \}$. One can imagine it as an indicator whether the edges were processed already along the transformation of the realization $X$ into realization $Y.$ (Green edges are ready for processing while red edges are processed already.)

For $1 \le i<j\le m$ denote $\green_f[e_i,e_j]$  the (not necessarily consecutive) subsequence of $\green_f$ edges from the sequence $e_i\cdots  e_j$. (The notation $\green_f$ is a shorthand for $\green_f[e_1,e_m],$ and the notations $\red_f [e_i,e_j]$ are defined analogously.) We will maintain the following property along our algorithm:
\begin{description}
\item[$(\pounds)$] any maximal consecutive $\red$ subsequence in $\mb e$ forms a closed alternating walk.
\end{description}
Let $f$ be a coloration on the current alternating walk $\mb e$ satisfying property $(\pounds)$. Furthermore let
\begin{equation}\label{eq:min-j}
j=\min\{j:\exists i<j\ \ \green_f[e_i,e_j]  \text{ is a cycle}\},
\end{equation}
and let
\begin{equation}\label{eq:min-i}
i=\max\{i:\ \green_f[e_i,e_j] \text{ is a cycle}\}.
\end{equation}
Since $\mb e$ is not necessarily a closed walk therefore such $j$ does not always exist. However if $\mb e$ is a closed walk and $\green_f$ is not empty, then such $j$ exists. Indeed, if a closed walk is deleted from a bigger closed walk then the remnant is a closed walk. Furthermore it is clear that integer $j$ determines uniquely the green cycle ending at $e_j.$ However before this cycle (along the original circuit) there may be several red edges. Therefore there may be several different integer $i$ defining the same green cycle. One way to handle this fact is equation (\ref{eq:min-i}).

\medskip\noindent Now we are ready to introduce our main tool to control the decomposition of an alternating circuit into alternating cycles.

For that end we define the operator $\tope$ on the edges of  walk $\mb e$ and the current coloration $f$ (satisfying condition $(\pounds)$) as follows:
\begin{definition}\label{operator}
$\tope({\mb e},f)$ will be a triple $({\mb e'},f',\mc C')$, where
\begin{enumerate}[(i)]
\item $\mc C'$ is the alternating cycle in $\mb e$ defined by equalities  (\ref{eq:min-j}) and (\ref{eq:min-i}), so
$$
 \mc C'=\green_f[e_i,e_j];
$$
\item $f':\mb e'\to \{\green,\red\}$ is defined with
$$
\red_{f'}=\red_f\cup \mc C';
$$
\item $\mb e'$ is an alternating walk obtained by rearranging the edges from $\mb e$.
\end{enumerate}
\end{definition}

\medskip\noindent
If $j$ is undefined, then $\tope(\mb e,f)$ is undefined. Let us remark that the length of $\mc C'$ is even, because $(K,F\cup F')$ was a bipartite graph, so $\mc C'$ is an alternating cycle.

What is missing is the description of the new alternating walk $\mb e'.$ Next we do just that. (Let's recall that two sequences written next to each other denotes their concatenation.) Write
$$
[e_i,e_j]=\mb g_1  \mb r_1 \cdots   \mb r_{k-1} \mb g_{k},
$$
where
$$
\green_f[e_i,e_j]= \mb g_1 \mb  g_2 \cdots   \mb  g_{k-1} \mb   g_{k}
$$
and
$$
\red_f [e_i,e_j]=\mb r_1 \mb  r_2 \cdots \mb r_{k-1}.
$$
In words: $\mb g_i$s and $\mb r_i$s represent the maximal consecutive $\green_f$ and $\red_f$ subsequences.
Let
$$
i'=\left\{
\begin{array}{ll}
\min\limits_\ell \{\ell <i: [e_{\ell},e_{i-1}]  \text{ is } \red_f \}, & \text{if }f(e_{i-1})=\red_f, \\
i, & \text{otherwise.}
\end{array}
\right.
$$
Furthermore let
$$
j'=\left\{
\begin{array}{ll}
\max\limits_\ell \{\ell>j:\ [e_{j+1},e_{\ell}] \text{ is } \red_f \}, & \text{if }f(e_{j+1})=\red_f, \\
j, & \text{otherwise.}
     \end{array}
   \right.
$$
We define
$$
\mb r^-=\left\{
\begin{array}{ll}
[e_{i'},e_{i-1}],&\text{if } i'<i, \\
\emptyset, &\text{if } i'=i;
\end{array}
\right .
$$ 
and
$$ 
\mb r^+=\left\{
\begin{array}{ll}
[e_{j+1},e_{j'}],&\text{if } j'>j, \\
\emptyset, &\text{if } j'=j.
\end{array}
\right .
$$ 

\medskip\noindent Let
\begin{equation}\label{eq:order}
\mb e'=e_1\cdots e_{i'-i}\mb r^+  \inv{\mb g_{k}}  {\mb r_{k-1}}
\inv{\mb g_{k-1}} \cdots  {\mb r_{1}}  \inv {\mb g_{1}}\mb r^-  e_{j'+1}\cdots e_m.
\end{equation}
This last formula requires some explanation: the cycle $\mc C'$ consists of the $\green_f$ segments of $[e_i,e_j].$ All the $\red_f$ segments form alternating closed walks that were processed earlier. We may assume without loss of generality, that the very first edge $e_i$ belongs to $F,$ consequently the last edge $e_j$ belongs to $F'.$

When we finish the required swap operations exchanging the edges from $F$ into edges from $F'$ along cycle $\mc C'$ (and transferring the actual degree realization closer to realization $Y$), then listing the edges of $\mc C'$ in the same way as before would not produce an alternating closed walk anymore. To form an alternating walk again we must consider the edges of $\mc C'$ in the opposite order. This is done by the subsegments $\inv{\mb g_i}$s. Listing $\mc C'$ in opposite order must list the closed walks $\mb r_i$s also in opposite order (see (\ref{eq:order})), which in turns takes care automatically for keeping the alternating order of edges from $F$ and $F'$.

One can ask the reason to exchange $\mb r^-$ and $\mb r^+$ since this is not necessary to keep the walk alternating.  This reason lays in equation (\ref{eq:orig-e}).

\bigskip\noindent Now we are ready to describe the control mechanism  to govern the swap sequence to change the edges of the current realization belonging to $F$ into the edges belonging to $F'$ along the alternating closed walk $\mb e.$ For that end let  $\green_f:=[e_1,e_m]$, furthermore let $\mb e_0:=\mb e$ and $f_0:= f$. Now we define the sequence
$$
({\mb e_1},f_1,\mc C_1), ({\mb e_2},f_2,\mc C_2),\dots, ({\mb e_n}, f_n,\mc C_n)
$$
by the formula
$$
({\mb e_{\ell+1}},f_{\ell+1},\mc C_{\ell+1}):=\tope ({\mb e_{\ell}},f_{\ell})
$$
for $\ell=0,1, \ldots$. We stop when $\tope ({\mb e_{n}},f_{n})$ is undefined.
We also define the sequence
$$
(F_0, F'_0),(F_1, F'_1),\dots, (F_{n}, F'_{n})
$$
of partitions of $F\cup F'$ as follows:
\begin{enumerate}[(1)]
\item let $F_0:=F$ and $F_0':=F'$,
\item let $F_{i+1}:=F_i\cup(\mc C_{i+1}\setminus F_i)\setminus  (\mc C_{i+1}\cap F_i)$ and $F'_{i+1}:=(F\cup F')\setminus F_{i+1}$.
\end{enumerate}
We define $n(\mb e):=n$ and
\begin{equation}\label{eq:nTf}
\iter^i(\mb e):=\mb e_i, \quad f^i(\mb e):=f_i \quad \text{ for } 0\le i\le n(\mb e).
\end{equation}
\old{
furthermore
\begin{equation}
\iter(\mb e):=\mb e_{n(\mb e)}\quad \text{ and }\quad f(\mb e):=f_{n(\mb e)}.
\end{equation}
    }
It is easy to see, and we will  show formally in Lemma \ref{lm:e-property}, that if $\mb e$ is a circuit then $\mc C_1,\dots, \mc C_{\mb e(n)}$ will be a circuit decomposition of $\mb e$. Later we will use this decomposition to obtain our canonical path system.

Let us emphasis here that we do not have an operation $\iter$, and $\iter^{i+1}({\mb e})$ is not computable from $\iter^{i}({\mb e})$ without knowing exactly which circuits of the current alternating walk have been processed.

We will prove a series of observations. We start with some easy direct consequences of definitions (\ref{eq:min-j}), (\ref{eq:min-i}) and (\ref{eq:order}):
\begin{lemma}\label{lm:not-before-red}
During the algorithm, at any given iteration $\kappa$ we have:
\begin{enumerate}[{\rm (i)}]
\item in the current alternating walk $\mb e_{\kappa -1}$ the edge  $e_{j_{\kappa}}$ is after all $\red_{f_{\kappa-1}}$ edges;
\item for any red edge the size of the maximal red subsequence containing it cannot decrease;
\item the number of maximal red subsequences can be increased by at most one, but can be decreased by any reasonable number. \qed
\end{enumerate}

\end{lemma}

\begin{lemma} \label{lm:e-property}
For each $0\le {\nu}\le \mb e(n)$ we have:
\begin{enumerate}[{\rm (i)}]
\item maximal $\red_{f_{\nu}}$  intervals $[\mb e_{\nu}(k), \mb e_{\nu}(\ell)]$ in  $\mb e_{\nu}$ are circuits $($recall, $\mb e_{\nu}(d)$ is the $d$th edge along $\mb e_\nu )$;
\item the edge sequence  ${\mb e}_{\nu}$ is a walk which alternates between $F_{\nu}$ and $F'_{\nu}$;
\item $\verti{\mb e_{\nu}}1=\verti{\mb e}1$ and $\verti{\mb e_{\nu}}{m+1} =\verti{\mb e}{m+1}$ $($these are the very first and very last vertices in $\mb e)$;
\item $\green_{f_{\nu}}[{\mb e_{\nu}}]$ is a walk from $v_{\mb e}(1)$ to $v_{\mb e}(m+1)$ $($while only just a part of the edges of $\mb e_{\nu}$ are $\green$ they still provide an alternating walk between those vertices$).$
\item if $\mb e$ is circuit, then $f(\mb e)$ is the constant red function, i.e. we processed all edges, while  $\mc C_1,\dots, \mc C_n$ is an alternating cycle decomposition of $\mb e$.
\end{enumerate}
\end{lemma}
\begin{proof}
We prove the statements by induction on ${\nu}$. For ${\nu}=0$ the statements are trivial because $\green_{f_0}=[e_1,e_m]$. Consider now the inductive step ${\nu}-1 \to {\nu}$. Assume that
\begin{equation}\label{eq:inti}
\mb e_{\nu-1}=e_1\cdots e_{i'} \mb r^- {\mb g_{1}}  {\mb r_{1}} \cdots  {\mb g_{k-1}}  {\mb r_{k-1}}  {\mb g_{k}} \mb r^+  e_{j'}\cdots e_m
\end{equation}
and
$$
\mb e_{{\nu}}=e_1\cdots e_{i'} \mb r^+ \inv{\mb g_{k}}  {\mb r_{k-1}}  \inv{\mb g_{k-1}} \cdots  {\mb r_{1}}  \inv{\mb g_{1}} \mb r^- e_{j'}\cdots e_m.
$$
(Here it is important to recall, that when $\mb r^-$ and/or $\mb r^+$ is empty, then $i'=i-1$ and/or $j'=j+1.$ If some of these cases apply, then the corresponding remarks on $\mb r^-$ and $\mb r^+$ are void.)

\smallskip
\noindent (i) The intervals $\mb r_\ell$ are maximal red intervals, so by the inductive assumption they are circuits, i.e. $\bott{\mb r_\ell}=\topp{\mb r_\ell}$. Moreover, by the construction, the first vertex of $\mb g_1$ and the last vertex of $\mb g_k$ are the same: $\bott{\mb g_1}= \topp{\mb g_k}$, and $\mb g_\ell$ is a  path from $\bott{\mb g_\ell}$ to $\topp{\mb g_\ell}$. Since $\topp{\mb g_\ell}=\bott{\mb r_\ell}= \topp{\mb r_\ell}= \bott{\mb g_{\ell+1}}$, we have that
$$
\mb c\ =\ \inv{\mb g_{k}}  {\mb r_{k-1}}  \inv{\mb g_{k-1}}  \cdots{\mb r_{1}}  \inv{\mb g_{1}}
$$
is a $\red_{f_\nu}$ circuit.

To finish the proof of (i) there is only one remaining case: if the maximal $\red_{f_{\nu}}$ interval $[e_l,e_\ell]$ in $\mb e_{\nu}$ contains properly the interval  $[e_i,e_j].$ (This is the case when at least one of $\mb r^-$ and $\mb r^+$ are not empty.) Then both $[e_l,e_{i-1}]$ and $[e_{j+1},e_\ell]$ are maximal $\red_ {f_{\nu-1}}$  intervals in $\mb e_{\nu-1}$ so $[e_l,e_{i-1}]=\mb r^-$ and $[e_{j+1},e_\ell]=\mb r^+$, therefore $[e_l,e_{i-1}] = \mb r^-[e_i, e_j]\mb r^+$ is  the concatenation of  at most three circuits  (since $\mb r^-$ or $\mb r^+,$ but not both, can be empty) , so it  is also a circuit.

\smallskip\noindent (ii) The vertices $v_{\mb r^-}(1)= v_{\mb r^+}(1)= v_{\mb g_1}(1)$ are identical in $\mb e_{\nu-1}.$ Therefore edges $e_{i'}$ and $e_{j+1}$ belong to $F_{\nu -1}$ simultaneously and the same applies for the edges $\last{\nu-1}{\mb r^-}=e_{i-1}$ and $\last{\nu-1}{\mb r^+}=e_{j'}$. (Here the index in  $\last{\nu-1}{}$ refers to the order of the walk $\mb e_{\nu-1}$.)  So, since $e_1\cdots e_{i'-1} \mb r^-$ an alternating walk in $\mb e_{\nu-1}$ therefore the same applies for $e_1\cdots e_{i'-1} \mb r^+$ (and analogously for $ \mb r^- e_{j'+1}\cdots e_m$) in $\mb e_{\nu}.$ In other words it makes no difference in the behavior (relating to the sub-walk $[e_i,e_j]$) of the walks $[e_1,e_{i-1}]$ and $[e_{j+1},e_m]$ whether $\mb r^-$ and/or $\mb r^+$ is/are empty.

Furthermore we have
$$
\topp{e_1\cdots e_{i'-1} \mb r^+}=\bott{\mb g_1}=\bott{\mb c}=\topp{\mb g_k}=\topp{\mb c} =\bott{\mb r^- e_{j'+1}\cdots e_m},
$$
so  $\mb e_{\nu}$ is a walk.

Next we check whether ${\mb e}_{\nu}$ alternates between $F_{\nu}$ and $F'_{\nu}$. Since we have $F_{\nu} \cap \{e_0\cdots e_{i'-1} \cup \mb r^+\}= F_{\nu-1} \cap \{e_0\cdots  e_{i'-1} \cup \mb r^+\}$, the interval $e_0\cdots e_{i'-1}\mb r^+$ alternates between $F_{\nu}$ and $F'_{\nu}$ and the analogous statements holds for $\mb r^- e_{j'+1}\cdots e_m.$

We know that
$$
e_{i-1} \in F_{\nu-1}  \Leftrightarrow  e_i \in F'_{\nu-1} \Leftrightarrow e_j  \in F_{\nu-1},
$$
since $[e_i,e_j]$ is a circuit in $\mb e_{\nu-1}.$ We also have  that
$$
e_{j'}\in F_{\nu} \Leftrightarrow e_i \in  F_{\nu} \Leftrightarrow e_j \in F'_{\nu}.
$$
Since $e_j$ is the first edge of $\inv{\mb g_k}$, the path $e_0\cdots e_{i'-1} \mb r^+ \inv{\mb g_{k}}$ alternates  between $F_{\nu}$ and $ F'_{\nu}$.

Assume that $\mb r_{k-1}=e_p\cdots e_r$. Then
$$
e_p\in F_{\nu-1} \Leftrightarrow e_r\in F'_{\nu-1} \Leftrightarrow e_{r+1}\in F_{\nu-1}.
$$
Thus
$$
e_p\in F'_{\nu} \Leftrightarrow e_r\in F'_{\nu} \Leftrightarrow e_{r+1}\in F_{\nu}.
$$
Therefore the path $e_1\cdots e_{i'-1} \mb r^+ \inv{\mb g_{k}} \mb r_{k-1}$ alternates  between $F_{\nu}$ and $ F'_{\nu}$ because  $\last{\nu}{\inv{\mb g_k}}$ is $e_{r+1}$ and $\first{\nu}{\mb r_{k-1}}$ is $e_p$.

Repeating the arguments above we obtain that the whole path $\mb e_{\nu}$ alternates  between $F_{\nu}$ and $ F'_{\nu}$ which finishes the proof of (ii).

\smallskip\noindent (iii) Here everything is trivial - except if $i=1$ and/or $j=m.$ By symmetry, it is enough to study one of these, let say $j=m.$ Then the last segment of $\mb e_\nu$ is $\mb r^+ \inv{\mb g_{k}}  \cdots   \inv{\mb g_{1}} \mb r^-$ which is a circuit, so the current end point of $\mb e_\nu$ is the same as the original end point of $\mb e_{\nu -1}.$

\smallskip\noindent (iv) All maximal $\red_{f_\nu}$ intervals are circuits, therefore removing them one by one  from $\mb e_\nu$ does not destroy the connectivity in $\green_{f_\nu}$ from $b(\mb e_\nu)$ to $t(\mb e_\nu)$ (as far as there are green edges).

\smallskip\noindent (v) It follows immediately from (iii) and (iv): a non-empty green remainder is a circuit, so the process will not finish while there still exists some green remainder. Consequently $\nu < \mb e(n).$
\end{proof}

\begin{lemma} \label{lm:g_r}
{\rm (a)} For each $0\le \nu\le n$ and $1\le r<s\le m$, if $\mb e_{\nu}(r)$ is $\green _{f_{\nu}}$ and $\mb e_{\nu}(s)$ is $\red_{f_{\nu}}$, then $\bott{\mb e_{\nu}(r)}\not \in \mb e_{\nu}(s) $.\\
{\rm (b)} Furthermore if $\mb e_{\nu}(r')$ is also $\green _{f_{\nu}}$ where  $ r<r'<s,$ then $\bott{\mb e_{\nu}(r)}\ne \topp{\mb e_{\nu}(r')}$.
\end{lemma}
\begin{proof}
Assume on the contrary that the statement is not true. Consider a counterexample
where $\nu$ is minimal. Assume that
$$
\mb e_{{\nu}-1}=e_1\cdots e_{i'-1} \mb r^- \mb g_{1}  \mb r_{1} \cdots  \mb g_{k-1}  \mb
r_{k-1}  \mb g_k  \mb r^+ e_{j'+1}\cdots e_m
$$
and
$$
\mb e_{{\nu}}=e_1\cdots e_{i'-1} \mb r^+  \inv{\mb g_{k}}  {\mb r_{k-1}}  \inv{\mb g_{k-1}}
\cdots  {\mb r_{1}}  \inv{\mb g_{1}} \mb r^- e_{j'+1}\cdots e_m.
$$
The edge sequence ${\mb g_{1}}  {\mb r_{1}} \cdots  {\mb g_{k-1}}  {\mb r_{k-1}}
{\mb g_{k}}$  (in $\mb e_{\nu-1}$) is a circuit.

Since  $\mb e_{\nu}(r)$ is unprocessed in $\mb e_{\nu}$ therefore $\mb e_{\nu}(r)\in e_0\dots e_{i'-1} \cup e_{j'+1}\cdots e_m.$ Furthermore $\mb e_{\nu}(s)\in \inv{\mb g_{k}}  \cup  \inv{\mb g_{k-1}} \cup \cdots \cup  \inv{\mb g_{1}}$ otherwise its color would be the same under $f_{\nu-1}$ and $f_{\nu}$ therefore $\nu$ would not be a minimal counterexample. But then the property $r< s$ infers that $\mb e_{\nu}(r)\in e_0\dots e_{i'-1}$.

Moreover $b(\mb e_{\nu -1}(r) )=b(\mb e_{\nu}(r)) \ne b(e_i)=t(e_j).$ Indeed, if $\mb r^-$  is not empty, then  $\mb e_{\nu -1}(r)$ and $e_{i-1}$ already would form a forbidden configuration in $\mb e_{\nu -1}$, a contradiction (the other case is similar). If both  $\mb r^-$  and  $\mb r^+$  are empty, then $[\mb e_{\nu-1}(r),\mb e_{\nu-1}(i-1) ]$ would be a circuit and it would contain a $\green_{f_{\nu}}$ cycle, a contradiction to the definition of $\mc C_{\nu}$ (in $\mb e_{\nu-1}$).

Therefore $b(e_{{\nu}}(r))$ must be an inner vertex of the cycle $e_i\dots e_j$. Now if this vertex is not the last vertex of a $\inv{\mb g_\ell}$, that is we have  $b(e_{{\nu}}(r))= b(e_{\nu}(s))$ then $\green_{f_{{\nu}-1}}[\mb e_{\nu-1}(r),\mb e_{\nu-1}(s)]$ would be a circuit, containing a  $\green_{f_{{\nu}-1}}$ cycle with smaller maximal element than $e_j$, which contradicts to the definition of $\mc C_{\nu}$ in $\mb e_{\nu -1}.$ Finally, if  this vertex is the last vertex of a $\inv{\mb g_\ell}$ then it is also the first vertex of $\mb r_{\ell-1}$ therefore edges $(e_{{\nu}}(r))$ and $\first{}{\mb  r_{\ell-1}}$ would form already in $\mb e_{\nu-1}$ the forbidden configuration of the statement, contradicting the minimality of $\nu.$

The proof of (b) uses a similar argument.
\end{proof}

\begin{lemma}\label{lm:eTre}
Assume that
\begin{equation}\label{eq:Tre}
 \iter^r(\mb e)=\mb g_1\mb r_1\mb g_2\mb r_2\dots\mb r_k\mb g_{k+1},
\end{equation}
where the first and/or the last green subsequence can be empty. Then
\begin{equation}\label{eq:orig-e}
 \mb e=\mb g_1\inv{\mb r_1}\mb g_2\inv{\mb r_2}\dots\inv {\mb r_k}\mb g_{k+1},
\end{equation}
so we obtain back the original edge sequence $\mb e$.
\end{lemma}
\noindent
It is important to understand that here we do not have any realization in the background (and no alternation is considered on the edges), we consider only the order of the edges. The operations above are nothing else, just turning back all maximal $\red_{f^r(\mb e)}$ intervals in $\iter^r(\mb e)$.
\begin{proof}
We apply mathematical induction on $r$. For $r=0$ the statement is trivial because $\iter^0(\mb e)=\mb e=\mb g_1$.

Now we assume that the statement is true for  $(r-1)$ and we are going to prove it for  $r$. For that end assume that
\begin{equation}\label{eq:van-r+}
\iter^{r-1}(\mb e)=\mb g_1\mb r_1\cdots \underbrace{\mb r^- \mb g_{t}\mb r_t \cdots\mb r_{u-1} \mb g_u \mb r^+} \mb g_{u+1} \cdots \cdots\mb r_k\mb g_{k+1}.
\end{equation}
where the formulas \ref{eq:min-j} and \ref{eq:min-i} select the intervals $\mb r^- \mb g_t\mb r_t \cdots\mb r_{u-1} \mb g_u \mb r^+$ to process (where  $\mb r^-$ and/or $\mb r^+$ can be  empty).

To compute $\iter^{r}(\mb e)$ we should check if $\mb r^-$ and $\mb r^+$ are empty or not. Altogether there are four cases to investigate, however the properties of one end of the sequence of $\mc C_r$ does not influence the other end, therefore it is enough to consider one  ``generic case'', say,  when $\mb r^-$ is empty but $\mb r^+$ is not empty. Then
\begin{equation}
 \iter^{r}(\mb e)=\mb g_1\mb r_1\cdots \underbrace{\mb r ^+ \inv{\mb g_u}\mb r_{u-1}\inv{\mb g_{u-1}}\cdots\mb r_{t} \inv{\mb g_t}}_{\text{red in $f^r(\mb e)$}} \mb g_{u+1}  \cdots \cdots\mb r_k\mb g_{k+1}.
\end{equation}
Now ${\mb r^+ \inv{\mb g_u}\mb r_{u-1}\inv{\mb g_{u-1}}\cdots\mb r_{t} \inv{\mb g_t}}$  is
a maximal red interval in $f^r(\mb e).$ When we ``turn back'' the $f^r(\mb e)$-red maximal intervals in $\iter^r(\mb e)$ we get:
\begin{eqnarray}
\mb g_1\inv{\mb r_1}\cdots {\mb g_{t-1}} \Big(\inv{\mb r^+ \inv{\mb g_u}\mb r_{u-1}\cdots\mb r_{t} \inv{\mb g_t}}\Big) \mb g_{u+1}\inv{\mb r_{u+1}} \cdots \cdots\inv{\mb r_k}\mb g_{k+1} &=& \nonumber \\
\mb g_1\inv{\mb r_1}\cdots {\mb g_{t-1}}\Big({\mb g_t} \inv{\mb r_{t}}\cdots {\mb g_{u-1}}
\inv{\mb r_{u-1}} \mb g_u\inv{\mb r^+}\Big)\mb g_{u+1} \cdots \cdots\inv{\mb r_k}\mb g_{k+1} &=&  \nonumber \\
\mb g_1\inv{\mb r_1}\cdots \mb g_{t-1} \mb g_{t} \inv{\mb r_{t}}\cdots {\mb g_{u-1}} \inv{\mb r_{u-1}}
\mb g_{u}\inv{\mb r^+} \mb g_{u+1}\cdots \cdots\mb r_k\mb g_{k+1}& = &\!\!  \mb e \label{eq:ize2}
\end{eqnarray}
where (\ref{eq:ize2}) is just the inductive assumption.
\end{proof}

\begin{lemma}\label{lm:Tte-new}
Assume that {\rm (\ref{eq:Tre})} holds, and define  $n_\ell:=n(\mb r_\ell)$ for $1\le \ell\le k$. Furthermore let $t_{\ell}:= r_1+ \cdots  + r_{\ell}$ for all $1\le \ell\le k$. Then
$$
\iter^{t_{\ell}} (\mb e')=
\mb g_1\inv{\mb r_1}\dots  \mb g_{\ell}\inv{\mb r_{\ell}} \mb {g_{\ell+1}}\mb r_{\ell+1}\cdots{\mb r_k}\mb g_{k+1}.
$$
\end{lemma}
\begin{remark}
It is important to emphases  that there is no reason that the algorithm running on $\mb e'$ would provide the same cycle decompositions of circuits $\mb r_i$ as the the same algorithm, running on the original $\mb e$ would do. As a matter of fact one can construct example where this is not the case.
\end{remark}
\begin{proof}[Proof of the Lemma \ref{lm:Tte-new}]
We apply induction on $\ell$. For $\ell=0$ there is no processed edge in $\mb e'$, nothing to prove. So assume that $\ell \ge 1$ and we know the statement for $\ell -1$.  For $v=0,\ldots, r_{\ell}$ let $\tau _v :=  t_{\ell-1} + v.$  We are going to show that
\begin{equation}\label{eq:partly}
\iter^{\tau_v}(\mb e')=
\mb g_1\inv{\mb r_1}\cdots \inv{\mb r_{\ell-1}} \mb g_{\ell}\iter^{v}(\mb r_\ell)\mb {g_{\ell+1}}\mb r_{\ell+1}\cdots \mb g_{k+1}.
\end{equation}
In words:  iterations $t_{\ell-1}\! +\!1, \dots, t_{\ell}$ of our algorithm work on $\mb r_{\ell}$ and completely process it, furthermore at each iteration we have
\begin{equation}\label{eq:partly1}
f^{\tau_v}(\mb e') \big | _{ \displaystyle \mb r_\ell} = f^v(\mb r_\ell).
\end{equation}
We prove it with induction on $v.$ When $v=0$ then we have nothing to prove, since  case $\tau_0$ coincides with $t_{\ell-1}.$ Assume now that (\ref{eq:partly}) holds for $\tau_{v-1}$ and prove it for $\tau_v.$

Compute $\iter^{\tau_v}(\mb e'_{\tau_v-1})$. By Lemma \ref{lm:not-before-red} (i) the current $e_{j_{\tau_v}}$ is after all $\red_{f_{\tau_v-1}}$ edges. However it is within $[r_{\ell}]_{\mb e'}$ since the original execution of our algorithm producing $\mb e'$ fully processed the closed walk $[r_{\ell}]_{\mb e'}$ while in $\mb e'_{\tau-1}$ it is not  achieved yet: there exists at least one not processed cycle. Finally, for the same reason, $e_{i_{\tau_v}}$ also should be in $[r_{\ell}]_{\mb e'}.$ So we know that  $[e_i,e_j]_{\tau}$ is a consecutive subset of of $\iter^{\tau_{v-1}}(\mb r_{\ell})$. (It is clearly not necessarily a subsequence!)

By the inductive hypothesis and (\ref{eq:partly1}) for $\tau_v-1$ we have:
$$
\green_{f^{v-1}(\mb r_\ell)}[e_i,e_j]= \green_{f^{\tau_v-1}(\mb e')}[e_i,e_j],
$$
therefore cycle $[\mc C_{\tau}]_{\mb e'_{\tau_v-1}}$ coincides with $[\mc C_{v}]_{\mb r_{v-1}}$. This proves  (\ref{eq:partly1}) for $v$ which, in turns, proves ( \ref{eq:partly}) for $\tau_v.$
\end{proof}
\medskip\noindent
Now we are ready to formalize the center piece of our control mechanism to govern the construction of the required multicommodity flow (or, in other words, the swap sequences between different realizations). With the previous definitions one can quantify the size of a parameter set to follow the current status of the cycles in the decomposition of the alternating circuit $\mb e.$ It clearly can be exponentially big, so this cannot prove fast mixing time.

However, we do not need to know the status of those cycles. What we really have to know is the original  walk $\mb e.$ And, surprisingly enough, we can determine it with high probability. More precisely the following property holds:
\begin{theorem}\label{tm:TiTi}
If $\mb e$ is a circuit, and $0\le s\le n(\mb e) $, then
\begin{equation}\label{eq:reconst}
   \iter^{s}(\iter^r(\mb e))=\mb e.
\end{equation}
for some $0\le s\le n(\iter^r(\mb e))$.
\end{theorem}
\begin{proof}
Write $\mb e'=\iter^r(\mb e)$ and assume {\rm (\ref{eq:Tre})} that is
$$
 \mb e'=\mb g_1\mb r_1\mb g_2\mb r_2\dots\mb r_k\mb g_{k+1}.
$$
The application of Lemma \ref{lm:Tte-new} for $\ell=k$ proves the statement.
\end{proof}
What this statement says is the following. Assume that we performed a certain amount of swaps  along the cycle decomposition of the original alternating circuit (using our decomposition algorithm) and we have the alternating circuit $\iter^r(\mb e)$ in hands. Then, if we consider this alternating circuit as a totally fresh one and we use our decomposition algorithm, furthermore we perform our swap operations along this decomposition, then this procedure will process the red $\mb r_{\ell}$ subsequences one by one. But our problem here is that we do not know - yet -  when this procedure processes  fully all necessary $\mb r_k$s. In other words: when we should halt the algorithm.

However knowing the number of processed edges in the fully processed circuits of $\mb e'$ fully solves this problem, since we can use this parameter to halt our algorithm on $\mb e'.$
And the size of the set of the possible numbers is simply linear. This set together with the polynomial running time of the algorithm named in (\ref{eq:reconst}) provides a polynomial mean to determine $\mb e$ with its alternations.

One can ask the reason why this newly developed method is so effective. In the attempted approach  described shortly at the beginning of Subsection \ref{ss:rec} we tried to deal with all possible cycle decompositions of the circuits (this is consist of all cycles and all their order). In the chosen algorithm, the analysis of it requires to consider only a quadratic number of possible cycle decompositions.

\subsection{Construction}

If $X,Y\in V(\G)$ let $E(X\bigtriangleup Y)$ be the symmetric difference of the edge sets $E(X)$ and $E(Y)$, set $E(X-Y)=E(X)\setminus E(Y)$, and $E(Y-X)=E(Y)\setminus E(X)$.

\medskip\noindent
Before we describe the construction of our multicommodity flow we need some further definitions:

\begin{definition}
For $T\in V(\G) $ let $M_T$ be the bipartite  $k\times l$ adjacency matrix of $T$. For
$X,Y,Z\in V(\G)$ write $\widehat  M(X+Y-Z) =M_X+ M_Y-M_Z$. (As we will see in the proof of the Key Lemma, these $k\times l$ matrices essentially {\em encode} the paths from $X$ to $Y$ along $Z.$)
\end{definition}

If $M$ and $M'$ are  $m\times m'$ matrices then let $\mathfrak{d}(M,M')$ be the
number of non-zero elements in $M-M'$ (the well-known {\em Hamming distance}).

\newcommand{\less}{\preceq'}
\newcommand{\lesss}{\preceq^*}

\medskip\noindent {\bf Outline of the construction of the path system.}
Fix a total order $\preceq$ on $U \times V.$ This will induce a total order $\less$ on all subsets of that product (namely we take the induced lexicographic order), in particular also on circuits in $[U,V]$. This will also induce a total order $\lesss$ on all sets of circuits in  $[U,V]$ (we can take again the induced lexicographic order).

\medskip\noindent For each $X \ne Y \in¸ V (\G)$ do the following.

\begin{enumerate}[(A)]
\item Let $S_{X,Y}=\S(E(X-Y), E(Y-X))$. (This notation was introduced at (\ref{eq:matching}).) To each  $s\in \S(E(X-Y), E(Y-X))$
consider the {\em unordered} alternating  circuit decomposition $\mc C_s$ of
$(E(X-Y), E(Y-X))$. (This is described in Lemma \ref{lm:decomp}.)
\item Order $\mc C_s$ using $\less$ to obtain the  {\em ordered} alternating
circuit decomposition
    \begin{displaymath}
            W^s_1,W^s_2\dots, W^s_{k_s}
    \end{displaymath}
    of $(E(X-Y), E(Y-X))$.
\item Every  $W^s_i$ is an alternating circuit in the bipartite graph $(U\cup V, E(X-Y)\cup E(Y-X))$. Consider the   enumeration $e_1\dots e_m$ of $W^s_i$, where ${e_1}$ is the $\prec'$-minimal edge in $W^s_i$,  and $e_2$ is the smaller edge for $\prec' $ among its two neighboring edges, while $e_m$ is the bigger. (This fixes uniquely the walk which traverses this circuit.) Now we can apply the method of Subsection \ref{ss:rec} to determine the cycle decomposition of $W^s_i$ for $1\le i\le k_s$:
    \begin{displaymath}
            C^{s,i}_1,C^{s,i}_2,\dots, C^{s,i}_{\ell_{s,i}}.
    \end{displaymath}
    Actually, we obtain  cycle $C^{s,i}_j$ as a sequence of edges. We keep this order to process  $C^{s,i}_j$ further in (F).
\item Let
        \begin{displaymath}
            C_1, C_2, \dots, C_{m_s}.
        \end{displaymath}
     be the short hand notation for the (alternating) cycle decomposition
        \begin{displaymath}
            C^{s,1}_1,C^{s,1}_2,\dots,   C^{s,1}_{\ell_{s,1}},C^{s,2}_1,C^{s,2}_2,\dots,
            C^{s,2}_{\ell_{s,2}},\dots, C^{s,k_s}_1,C^{s,k_s}_2,\dots, C^{s,k_s}_{\ell_{s,k_s}}
        \end{displaymath}
     of $E(X\bigtriangleup Y)$. We will call it a {\em canonical} cycle decomposition.
\item For each cycle $C$ in this decomposition we inherit an enumeration of that cycle (see (C)), which also determines a direction on the cycle. So for  $a,b\in C$ we can define  $[a,b]_C$ as the walk from $a$ to $b$ in $C$ according to this fixed direction.
\end{enumerate}
The following observation plays a crucial role in our method:
\begin{observ}\label{keyobserv}
The function $s$ itself determines this canonical decomposition, and also the direction of the cycles in the decomposition. So we do not need to know $E(X-Y)$ and $E(Y-X)$ to compute the ${C^{s,i}_j}$, or even $[a,b]_{C^{s,i}_j}$ from $s$.
\end{observ}

\begin{enumerate}[(F)]
\item Let  $\varUpsilon(X,Y,s)$ be a path of realizations
    \begin{equation}\label{eq:cycles}
         X=\mathrm{G}_0,\mathrm{G}_1, \dots,\mathrm{G}_{n_1}, \mathrm{G}_{n_1+1},\dots, \mathrm{G}_{n_2},\dots, \mathrm{G}_{n_{m_s}}=Y
    \end{equation}
    in $\G$ from $X$ to $Y$ such that
        \begin{enumerate}
            \item $n_{m_s}\le c\cdot n^2$,
            \item $E(\mathrm{G}_{n_i})=\big(E(\mathrm{G}_{n_{i-1}})\cup (E(Y)\cap E(C_i)\big) \setminus \big ( E(X)\cap E(C_i)\big )$ for $i = 1, 2, . . . ,m_s$,
            \item  if for $i< m_s$ we denote  the first vertex of the cycle $C_{i+1}$ in the order inherited from the construction by $a_{{i+1}}$, then for each $n_i\le j\le n_{i+1}$  there is is a vertex $b_j$ in $C_{i+1}$ such that
                \begin{displaymath}
                   \big  |E(\mathrm{G}_j)\bigtriangleup F \big | \le \Omega_1,
                \end{displaymath}
                where
                \begin{displaymath}
                    F=\big(E(\mathrm{G}_{n_i})\cup \big([a_{i+1},b_j]_{C_{i+1}}\cap
                    E(Y)\big)\setminus \big([a_{i+1},b_{j}]_{C_{i+1}}\cap
                    E(X)\big)\big),
                \end{displaymath}
            \item for each $j=1, 2, \ldots, n_{m_s}$ there is $T\in V(\G)$ such that
                \begin{displaymath}
                    \mathfrak{d}\left (\widehat M(X+Y-\mathrm{G}_j),M_{T}\right ) \le \Omega_2,
                \end{displaymath}
        \end{enumerate}
    where $c, \Omega_1$  and $\Omega_2$ are fixed ``small'' natural numbers. (Recall here, that by the definition of the Markov chain, in this path each graph $G_{\ell+1}$ is constructed from the previous one $G_\ell$ by a valid swap operation.)
\end{enumerate}

\begin{kobserv}\label{lm:keyob}
Let $X\ne Y \in \G$. If we can assign paths
\begin{displaymath}
\Bigl\langle\varUpsilon(X,Y,s):s\in \S(E(X-Y), E(Y-X)), X,Y\in V(\G) \Bigr\rangle
\end{displaymath}
according to  {\rm (A)-(F)} then (\ref{eq:multiflow_simple2}) holds and so
our Markov chain is rapidly mixing.
\end{kobserv}

\medskip
{\noindent {\bf Proof of the Key Lemma}:\ }
Fix $Z\in V(\mbb G)$. We need to prove (\ref{eq:multiflow_simple2}):
\begin{displaymath}
\sum_{\genfrac{}{}{0pt}{}{\scriptstyle X,Y \in     V(\mathbb{G})}{ }}
\frac{\left | \big \{s\in S_{X,Y}: Z\in \varUpsilon(X,Y,s) \big \} \right |
}{|S_{X,Y}|}
\le poly(n)\cdot N.
\end{displaymath}
\noindent  Let
\begin{displaymath}
\mathfrak{M}=\Bigl\{\widehat M(X+Y-Z):  Z\in \varUpsilon(X,Y,s)\text{ for some
} X,Y\in V(\G) \text{ and }\ s\in \S(\sDE)\Bigr \}.
\end{displaymath}
By  (F)(d) for each $\widehat M=\widehat M(X+Y-Z)\in \mathfrak{M}$
there is $T\in V(\G)$ such that  $\mathfrak{d}(\widehat M({X +Y-Z}), M_T)\le
\Omega_2$,
i.e. there are at most $\Omega_2$ positions  where $M_T$ and $\widehat M({X+Y
-Z})$  are different, so we have at most $ (n^2)^ {\Omega_2}= n^ {2 \Omega_2}$
difference sets.
Furthermore every entry of $\widehat M(X+Y-Z)$ lies in the set $\{-1, 0, 1,
2\}$,
so a fixed difference set we have most $3^{\Omega_2}$ possibilities. So
\begin{displaymath}
| \mathfrak{M}|\le |V(\G)|\cdot n^ {2 \Omega_2} \cdot 3^{\Omega_2}
\le \mathrm{poly}(n) \cdot |V(\G)|=\mathrm{poly}(n)\cdot N.
\end{displaymath}
For $\widehat M\in \mf M$ let
\begin{equation}
\mf X(Z, \widehat M)=\left \{(X,Y,s):s\in \S(X,Y),\ Z\in \varUpsilon(X,Y,s),\
\widehat M(X+Y-Z)=\widehat M \right \}.
\end{equation}
Since $|\mf M|\le \mathrm{poly}(n)\cdot N$, if we can prove that
\begin{equation}
\label{eq:enough1}
\sum_{\scriptstyle (X,Y,s) \in  \mf X(Z, \widehat M) }
\frac{1}{|S_{X,Y}|} \le poly(n)
\end{equation}
for all $M\in \mf M$, then  (\ref{eq:multiflow_simple2}) holds.

\bigskip\noindent
To verify (\ref{eq:enough1}) fix $\widehat M\in \mf M$.
Let  $(X,Y,s)\in \mathfrak{X} (Z,\widehat M)$ be arbitrary.
Since $M_Z+\widehat M=M_X+M_Y$,
we can compute
$$\Delta=E(X\bigtriangleup Y)$$
from $Z$ and $\widehat M$. Denote  by $(2d_1,\ldots,2d_h)$ the degree sequence
of $E(X\bigtriangleup Y)$. Put
\begin{displaymath}
t_\Delta=\prod_1^h (d_i!).
\end{displaymath}
Clearly
\begin{displaymath}
t_\Delta=\left |S_{X,Y}\right |,
\end{displaymath}
and so
\begin{displaymath}
\sum_{\scriptstyle (X,Y,s) \in  \mf X(Z, \widehat M) } \frac{1}{|S_{X,Y}|}=
\sum_{\scriptstyle (X,Y,s) \in  \mf X(Z, \widehat M) }
\frac{1}{t_\Delta}=\frac{\left |\mf X(Z, \widehat M)\right |}{t_\Delta}.
\end{displaymath}
Thus to prove (\ref{eq:enough1}) we need to show that
\begin{equation}\label{eq:enough2}
\left |\mf X(Z, \widehat M) \right |\le poly(n)\cdot t_{\Delta}.
\end{equation}
Let
\begin{equation}
\mc S=\left \{s: \mbox{ for some } (X,Y) \mbox{ we have } (X,Y,s)\in \mf X(Z, \widehat M)\right \}.
\end{equation}
To get (\ref{eq:enough2}) it is enough to show the following statement.
\begin{lemma}\label{lm:a-and-b}
For each possible $Z$ and the corresponding set $S$ we have:
\begin{enumerate}[{\rm (a)}]
\item $|\mc S|\le poly(n)\cdot t_{\Delta}$,
\item for each $s\in \mc S$ we have
\begin{equation}
\left |\left \{(X,Y):(X,Y,s)\in \mf X(Z, \widehat M)\right \}\right |\le
poly(n).
\end{equation}
\end{enumerate}
\end{lemma}
\noindent
To prove this lemma fix  $(X,Y,s)\in \mathfrak{X}(Z,\widehat M)$. We should recall the construction of the path $\varUpsilon(X,Y,s)$ which can be demonstrated as:
\begin{equation}
\overbrace{\ G^{k,\ell}_1,\ \cdots, \ G^{k,\ell}_m,\ \cdots ,\ G^{k,\ell}_{m_{k,\ell}}}^{
\overbrace{C^k_1,\ \quad \cdots\quad , C^k_\ell,\ \quad \cdots\quad , C^k_{\ell_k}}^
{\text{\normalsize $W_1,\ \quad \cdots\quad \cdots \quad W_k, \ \quad \cdots\quad \cdots \quad  W_{k_s} $}}
}
\end{equation}
where
\begin{enumerate}[(1)]
\item  we consider  first the circuit decomposition  of $(E(X-Y), E(Y-X))$ determined by $s$:
    $$
            W_1,W_2\dots, W_{k_s};
    $$
\item then, using the method of subsection \ref{ss:rec} for each $1\le k\le k_s$ we define an alternating circuit decomposition of $W_k$:
    $$
            C^{k}_1,C^{k}_2,\dots, C^{k}_{\ell_{k}};
    $$
\item  then in (F) for each $1\le k\le k_s$ and $1\le \ell \le {\ell_k}$ we define a sequence of elements of $\G$:
    $$
        G^{k,\ell}_1,\dots , G^{k,\ell}_m,\dots G^{k,\ell}_{m_{k,\ell}},
    $$
    such that
    \begin{multline}\label{eq:egk1}
        E\left (G^{k,\ell}_1\right )=\left [ E(Y)\cap \Big (\bigcup_{k'<k}E(W_{k'})\cup \bigcup_{\ell'<\ell} E\left  (C^k_{\ell'}\right ) \Big ) \right ] \ \bigcup \\
        \left [ E(X)\cap \Big (\bigcup_{k'>k}E(W_{k'})\cup \bigcup_{\ell'\ge \ell}E\left (C^k_{\ell'} \right )\Big ) \right ],
    \end{multline}
    and $G^{k,\ell}_{m_{k,\ell}}=G^{k,\ell+1}_1$ if $\ell<\ell_k$, $G^{k,\ell_k}_ {m_{k,\ell_k}}=G^{k+1,1}_1$ if $k<k_s$, and $G^{k_s,\ell_l}_{m_{k_s, \ell_k}}=Y$ (the equation \ref{eq:egk1} is just a reformulation of (F)(b));
\item finally  $\varUpsilon(X,Y,s)$ is the path
    \begin{equation}
        X=G^{1,1}_1,G^{1,1}_2, \dots, G^{k,\ell}_m ,\dots,  G^{k,\ell_k}_{m_{k,\ell_k}}=Y
    \end{equation}
    in $\G$ from $X$ to $Y$ (see (\ref{eq:cycles})).
\end{enumerate}
\noindent Fix $k$, $\ell$, $m$  such that $Z=G^{k,\ell}_m$, which means that  we are processing the $\ell$th cycle from the $k$th circuit.

By (F)(c) there are two vertices $a$ and $b$ in $C^{k}_{\ell}$ such that
\begin{equation}\label{eq:majdnemF}
 \big  |E(G^{k,\ell}_m)\bigtriangleup F \big | \le \Omega_1,
\end{equation}
 where
\begin{equation}\label{eq:F}
F=\left (E\left ({G^{k,\ell}_1}\right )\cup \big([a,b]_{C^{k}_{\ell}}\cap E(Y)\big)\setminus \big([a,b]_{C^{k}_{\ell}}\cap E(X)\big)\right ).
\end{equation}

\noindent To prove Lemma \ref{lm:a-and-b} (b) we show that
\begin{itemize}
\item[($\dag$)] \em there is a function $\Psi$ and a parameter set $\B$ such that $\B$ has $\mathrm{poly} (n)$ elements, and for each $(X,Y,s)\in \mathfrak{X}(Z,\widehat M)$  there is $B\in \B$ such that
    \begin{equation}\label{eq:teljesx}
        \Psi\left ( Z,\widehat M(X+Y-Z),s,B \right)= (X,Y).
    \end{equation}
\end{itemize}
Recall that $Z=G^{k,\ell}_m$ so we have $E\left (G^{k,\ell}_m\right )$. If we choose the parameter $B$ as the quadruple $\Bigl (i,a,b, $ $ E\left (G^{k,\ell}_m\right )   \Delta F \Bigr )$, then  using this parameter we can compute $F=E(G^{k,\ell}_m) \Delta \big (E(G^{k,\ell}_m)   \Delta F \big )$, and so
\begin{multline}\label{eq:exminusey}
E(X)\setminus  E(Y)=  \left ([a,b]_{C^{k}_{\ell}}\setminus  F\right )\bigcup \left ([b,a]_{C^{k}_{\ell}}\cap F\right )  \bigcup \\
\left [ \Big (\bigcup_{k'<k}E(W_{k'})\cup \bigcup_{\ell'<\ell} E\left  (C^k_{\ell'}\right )\Big )   \setminus F  \right ]\bigcup
\left [F \cap \Big (\bigcup_{k'>k}E(W_{k'})\cup \bigcup_{\ell'\ge \ell}E\left (C^k_{\ell'} \right )\Big ) \right ] .
\end{multline}
Since $i\le n^2$, $a,b\le n$ and $\big (E(X)\setminus E(Y) \big )\bigtriangleup F$ is an at most $\Omega_1$ element subset of $[U\cup V]^2$, the size of the parameter set is polynomial:
$$
|\mathbb B|\le n^2\cdot n\cdot n \cdot (n^2)^{\Omega_1}.
$$
Since $Z$ and $\widehat M(X+Y-Z)$  determine $E(X)\cap E(Y)$, we can compute $E(X)$. Similarly we can compute $E(Y)$. So we verified ($\dag$), and so Lemma \ref{lm:a-and-b} (b) holds.

\medskip\noindent
Now we turn to prove Lemma \ref{lm:a-and-b} (a).  We will do it in steps (a1) -- (a3).
\begin{enumerate}[{\bf ({a}1)}]
\item {\em Each function $s\in S,$ which corresponds to the  circuit decomposition
    $$
        W_1,\dots W_{k-1}, W_{k},  W_{k+1} \dots, W_{k_s}
    $$
    {\rm (}see Lemma \ref{lm:decomp}{\rm )}, is computable - using a small parameter set - from function $s'$ and the corresponding  circuit decomposition
    $$
        W_1,\dots W_{k-1},\iter^{\ell-1}( W_{k}),   W_{k+1} \dots, W_{k_s}.
   $$
   }
\end{enumerate}
Indeed,   by Theorem \ref{tm:TiTi}, $\iter^t (\iter^{\ell-1}( W_{k}))=W_{k}$ for some $t$.
So, as we described after Theorem \ref{tm:TiTi},  the  parameter pairs $k$ and the number of the processed edges in circuit $W_k$ together determine fully $\mb e$, and $k\le n^2$ and  the number of processed edges is also $\le n^2$. So {\bf (a1)} holds. \qed${}_{(a1)}$

\bigskip\noindent We need some preparation before we can formulate and prove (a2): Recall that Lemma \ref{lm:e-property} (ii) infers
$$
s'\in \mbb S\left (\Delta\cap E\left (G^{k,\ell}_1\right ),\Delta\setminus E\left(G^{k,\ell}_1\right ) \right ).
$$
The sequence
$$
 \mb e'=e_1e_2\dots e_m= \iter^{\ell-1}(W_k).
$$
is an alternating circuit in $G^{k,\ell}_1$. (All circuits of the decomposition with $k' < k$ are already fully processed. No circuit after $W_k$ is touched yet. So it is enough to consider only this.)  We use the notations $f=f^{\ell-1}(W_k)$ and
$$
\mb e'=\mb g_1\mb r_1\dots \mb g_u
$$
where $\mb g_s$ are maximal  $\green_f$, and $\mb r_o$ are maximal $\red_f$ intervals.
We know that the current cycle:
$$
C^k_\ell=\green_{f}[e_i, e_j]
$$
(for some $0\le i<j\le n$)  is undergoing a series of swaps operations which will exchange its edges between the realizations $X$ and $Y.$ When this swap sequence is completed then the processing of this cycle in  the cycle decomposition of circuit $W_k$ will be done, and the coloration of its edges will become $\red_{f^{\ell}}$. Now  the assumption (F)(c) about our swap sequence generation, applying for $G^{k,\ell}_m$, gives us an interval
$$
[a,b]_{C^k_\ell}=\green_{f}[e_i, e_{j'}],
$$
for some $i\le j'\le j.$

Assume that $e_i\in \mb g_o$ and $e_{j'}\in \mb g_t$. Write $\mb g_o=\mb  g_{o,0}\mb  g_{o,1}$, where $e_i$ is the first edge of $\mb  g_{o,1}$, and write $\mb g_t=\mb  g_{t,0}\mb  g_{t,1}$, where $e_{j'}$ is the last edge of $\mb  g_{t,0}$

Consider the sequence
$$
\mb e''= \mb g_1\mb r_1\dots \mb r_{o-1}\mb g_{o}\inv{\mb r_{o}}\mb g_{o+1}
\dots\inv{\mb r_{t-1}}\mb g_t \mb r_k\dots \mb g_u
$$
Now $\mb e''$ is in $\mbb (\Delta\cap F, \Delta\setminus F)$ because it is the concatenation of three $(\Delta\cap F, \Delta\setminus F)$-alternating paths. They are $\mb g_1\mb r_1\dots \mb r_{o-1}\mb g_{o,0}$ and  $\mb g_{o,1}\inv{\mb r_{o}}\mb g_{o+1} \dots\inv{\mb r_{r-1}}\mb g_{t,0}$ finally  $\mb g_{t,1}\mb r_{t+1}\dots \mb g_u$. However in general this walk is not necessarily alternating, because on the border of $\mb g_{t,0}$ and  $\mb g_{t,1}$ furthermore on the border of $\mb g_{t,0}$ and  $\mb g_{t,1}$ is not alternating anymore. (However, if one or both of the red circuits $\mb r^-$ and/or $\mb r^+$ exist then this problem will not occur there (see equation (\ref{eq:van-r+})).

\begin{enumerate}[{\bf ({a}2)}]
 \item {\em The walk $e'$ is computable {\rm (}using a small parameter set{\rm )} from $\mb e''$.}
\end{enumerate}
Indeed, $\mb e^+=\mb g_{o,1}\inv{\mb r_{o+1}}\dots \inv{\mb {r_{t-1}}}\mb g_{t,0}$ is computable from $\mb e''$ because it is a subsequence. Since
$$
[e_i,e_j]=\mb g_{o,1}\mb r_{o+1}\dots \mb {r_{t-1}}\mb g_{t,0}\mb g_{t,1}\dots \mb r_{t+1}\dots e_j
$$
is a circle, we can apply Theorem \ref{tm:TiTi} and  Lemma \ref{lm:Tte-new} to find some $v\le n^2$ such that $\iter^v (\mb e^+)=\mb g_{o,1}{\mb r_{o+1}}\dots {\mb {r_{t-1}}\mb g_{t,0}}$.
Thus  $[e_i, e_j]$ (in $\mb e'$) is computable from $\mb e''$. Since $\mb e'$ and $\mb e''$ agree    outside $[e_i, e_j]$, we proved {\bf (a2)}.

We turn our attention now to the third obstacle: till now we showed that knowledge of $[a,b]_{C^k_\ell}$ would determine fully $s$ from $s'.$ However we do not know exactly the sequence, since the assumption of (F) (c) allowed that $ \big  |E(\mathrm{G}_j)\bigtriangleup F \big | \le \Omega_1,$ so a small number of edges of the current realization is not on the alternating path determined by $s'.$ Next we will deal with this problem:

\begin{enumerate}[{\bf ({a}3)}]
\item {\em The sequence
    $$
        W_0,\dots, W_{\ell-1}, \mb e'', W_{\ell+1},\dots, W_k
    $$
    is  ``almost'' in $\mbb S\left ({\Delta}\cap E\left (G^{k,\ell}_m \right ),{\Delta}\setminus E\left (G^{k,\ell}_m\right ) \right )$, {\rm (}see formula \ref{eq:almost} below{\rm )} so it is computable from some element of $\mbb S\left ({\Delta}\cap E\left (G^{k,\ell}_m \right ),{\Delta}\setminus E\left (G^{k,\ell}_m\right ) \right )$ using a small parameter set.}
\end{enumerate}
Let
$$
\est=\left (E\left ({G^{k,\ell}_m}\right )\bigtriangleup F \right )\cup {C^k_{\ell}}(\bott{e_i})
\cup {C^k_{\ell}}(\topp{e_{j'}}).
$$
The last two expressions stand for the two edge pairs which are adjacent to the vertices of $\bott{e_i}$ and $\topp{e_{j'}}$ in the actual cycle ${C^k_{\ell}}.$ Since ${C^k_{\ell}}$ is a cycle indeed, we have  $|\est|\le   \Omega_1+4$ (due to (\ref{eq:majdnemF})).

For $w\in U\cup V$, let
$$
 t(w)=s''(w)\cap \left [\left ({\Delta}\cap E\left ({G^{k,\ell}_m}\right )\right )(w),\left ({\Delta}\setminus E\left ({G^{k,\ell}_m}\right )\right )(w)\right ],
$$
i.e $t(w)$ is those elements of $s''(w)$ which alternate between $\Delta\cap E\left ({G^{k,\ell}_m} \right )$ and   ${\Delta}\setminus E\left ({G^{k,\ell}_m}\right )$. Since $t(w)$ is a set of independent edges in
$$
H=\left [\left ({\Delta}\cap E({G^{k,\ell}_m})\right )(w),\left ({\Delta}\setminus E\left ({G^{k,\ell}_m} \right ) \right )(w)\right],
$$
we can find a perfect 1-factor extension $s^*(w)$ of $t(w)$  in the complete bipartite graph $H$. Then for this perfect matching we have
$$
s^* \in \mbb S\left ({\Delta}\cap E\left ( {G^{k,\ell}_m}\right ),{\Delta}\setminus E\left ({G^{k,\ell}_m}\right ) \right ).
$$
Next we show that the difference between $s''$ and $s^*$ is small, namely
\begin{equation}\label{eq:almost}
\sum_{w\in U\cup V} \left |s''(w)\bigtriangleup s^*(w) \right | \le 4\Omega_1 +16.
\end{equation}
For that end let $w\in U\cup V$ and $(e,e')\in s''(w)$ such that $e,e'\notin E^*$.
Then, by the definiton of $E^*$, we have
$$
e\in \Delta\cap G^{k,\ell}_m \Leftrightarrow e\in \Delta\cap F \quad\hbox{
and}\quad e'\in \Delta\cap G^{k,\ell}_m \Leftrightarrow e'\in \Delta\cap F.
$$
So since  $(e,e')$ alternates between  $\Delta\cap F$ and $\Delta\setminus F$,
$(e,e')$ also alternates between  $\Delta\cap G^{k,\ell}_m$ and $\Delta\setminus G^{k,\ell}_m$. So $(e,e')\in t(w)\subset s^*(w)$, and so $(e,e')\notin s''(w)\bigtriangleup s^*(w)$. Thus
\begin{multline}\label{eq:almost2}
\sum_{w\in U\cup V} \left |s''(w)\setminus s^*(w) \right | \le \\
\Big|\bigcup_{w\in U\cup V}\Big\{(w,e,e'):
(e,e')\in s''(w), e\in E^* \text{or } e'\in E^*\Big\}\Big|\le\\
2\cdot |E^*|\le  2\Omega_1 +8.
\end{multline}
Since $|s''(w)|=|s^*(w)|$, we have
$2\cdot |s''(w)\setminus s^*(w)|=|s''(w)\bigtriangleup s^*(w)|$,
so (\ref{eq:almost2}) gives  (\ref{eq:almost}).

Therefore $s^*$ together with a small parameter set which describes the symmetric
differences $s''(w)\bigtriangleup s^*(w)$ for $w\in U\cup V$
 determines completely $s''$,
thus  {\bf (a3)} is true as well.

\medskip\noindent
Putting together (a1)--(a3) we obtain
\begin{eqnarray*}
 |\mc S| & \le & poly_1(n)\cdot poly_2(n)\cdot poly_3(n)\cdot \left|\mbb S\big ({\Delta}\cap E(G^{k,\ell}_m),{\Delta}\setminus E(G^{k,\ell}_m) \big )\right| \\
& = & poly(n)\cdot t_\Delta.
\end{eqnarray*}
So  Lemma \ref{lm:a-and-b} (a) holds, which in turns completes the {\bf proof of the Key Lemma}. \qed

\bigskip\noindent We try to carry out the plan we just described. So:
\begin{itemize}
\item Fix $X \ne Y\in V(\G)$.
\item Pick $s\in \S(\sDE)$.
\item $s$ gives  an alternating  cycle decomposition
        \begin{equation}
            C_0,C_1,\dots, C_{\ell}
        \end{equation}
    of $E(X\bigtriangleup Y)$.
\end{itemize}
We want to define a path
\begin{equation}
  X=G_0, \dots,G_i,\dots, G_{m}=Y
\end{equation}
from $X$ into $Y$ in $\G$ - denoted by  $\varUpsilon(X,Y,s)$ -  such that
\begin{enumerate}[(i)]
\item the length of this path is $\le c\cdot n^2$ (where $c$ is a suitable
constant),
\item for some increasing indices $0<n_1<n_2<\dots n_\ell$ we have
$G_{n_i}=H_i$, where
    \begin{equation}
        E(H_i)=E(X)\bigtriangleup\left(\bigcup_{i'<i}E(C_{i'})\right).
    \end{equation}
\end{enumerate}
So we have certain ``fixed points'' of our path $\varUpsilon(X,Y,s)$, and this
observation reduces  our task to the following:
\begin{itemize}
\item for each $i<\ell$ construct the path
    \begin{equation}
        H_i=G'_0,G'_1,\dots, G'_{m'}=H_{i+1}
    \end{equation}
    between $G_{n_i}$ and $G_{n_{i+1}}$ such that $m' \le c\cdot |C_i|$ and
(F)(d)
holds, i.e. for each $j$ there is $K_j\in V(\G) $ such that $\mathfrak{d}\left
(\widehat
M(X,Y,G'_{j}),K_j \right )\le \Omega_2$.
\end{itemize}
{F}rom now on we work on that construction. To simplify the notation we write
$G=H_i$ and $G'=H_{i+1}$. We know that the symmetric difference of $G$ and $G'$
is just the cycle $C_i$. Now we are in the following situation:

\bigskip \goodbreak
\noindent{\bf Generic situation - construction of a path along a cycle}
\begin{enumerate}[(i)]
\item $X,Y,G,G'\in V(\G)$.
\item The symmetric difference of
$E(G)$ and $E(G')$ is a cycle $C$.
\item the symmetric differences $E(X\bigtriangleup G)$,
 $E(G\bigtriangleup G')$ and $E(G'\bigtriangleup Y)$
are pairwise disjoint.
\end{enumerate}
Construct a path
\begin{equation}
G=G_0,\dots, G_m=G'
\end{equation}
in the graph $\G$ of all realizations such that
\begin{enumerate}[(I)]
\item $m\le c\cdot|C|$, and the requirement of (F)(c) also holds,
\item for each $j$ there is $K_j\in V(\G)$ such that $\mathfrak{d}\left
(\widehat
M(X,Y,G_j), M_{K_j}\right )\le \Omega_2$.
\end{enumerate}
We will carry out this construction in the next sections. The burden of such a
construction is to meet requirement (II). In \cite{KTV97}  and in \cite{CDG07}
the regularity of the realizations was used.

\medskip

\noindent  The {\bf   ``friendly path method"}.

\noindent
In the next sections we describe a new general method based on the notion of
{\em friendly paths} (see Definition \ref{def:friendly}) to construct the paths
$\varUpsilon(X,Y,s)$.

The novelty of our friendly path method can be summarized as follows:
\begin{itemize}
\item if our bipartite  degree sequence is half-regular then the paths
$\varUpsilon(X,Y,s)$ satisfy the previous condition (II)
\item if our bipartite  degree sequence is arbitrary, then $\varUpsilon(X,Y,s)$
satisfies (II) provided the symmetric difference of $X$ and $Y$ is a cycle.
\end{itemize}
Originally we conjectured,  that our friendly path method  always produces paths
which satisfy (II). However we were unable to prove it, and now we think that
essentially new ideas are needed to prove the case of general bipartite degree
sequences.

\section{Multicommodity flow - along a cycle}\label{sec:bootleneck}
Let $X,Y$ and $Z$ be three realizations of a given bi-graphical degree sequence.  Assume that $E(X) \cap E(Y) \subset E(Z),$ furthermore $E(Z)\subseteq E(X)\cup E(Y)$. Then the realization $Z$ is an {\bf intermediate} realization between $X$ and $Y.$

In this section we describe the construction a path along an alternating cycle $C$. Here we have the intermediate realizations $G$ and $G'$ between $X$ to $Y$, and these two realizations differ only in this cycle $C,$ where $G\cap C = X\cap C$ and $G'\cap C = Y\cap C.$ At the beginning of this phase our canonical path is between  $X$ and $G$.   Along the process we extend it to reach realization $G'.$ Within the process all swaps will happen between vertices $V(C)$ of the cycle $C$ and the end of the process each chord will be at the same state as it was at the beginning, except the edges along the cycle, where the $X$-edges will be exchanged by the $Y$-edges.

In what follows we will imagine our cycles as convex polygons in the plane, and
we will denote by the vertices of any particular  cycle of $2\ell$ edges with
$u_1,v_1,u_2,$ $v_2,\ldots ,u_\ell, v_\ell$. The edges of the cycle are
$(u_1,v_1), (v_1,u_2), \ldots,(u_\ell,v_\ell),$  $(v_\ell,u_1)$ and they belong
alternately to $X$ and $Y.$ All the other (possible, but not necessarily
existing) edges among vertices of a particular cycle are the {\em chords}. (In
other words we will use the notion of chord if we want to emphasis that we do
not know whether the two vertices form an edge or not in the current graph.)  A
chord is a {\em shortest} one, if in one direction there are only two vertices
(that is three edges) of the cycle between its end points. The middle edge of
this three is the {\em root} of the chord.

W.l.o.g.  we may assume that $(u_1,v_1)$ is an edge in  $G$ while $(v_1,u_2)$ belongs to $G'.$  We are going to construct now a sequence of graphical realizations between $G$ and $G'$ such that any two consecutive elements in this sequence differ from each other in one swap operation. The general element of this sequence will be denoted by $Z.$

We have to control which graphs belong to this sequence. For that purpose we
assigned a matrix $\widehat M$ to each graph $Z.$ If $G$ is a vertex in $\G$ then  $M_G$
denotes the adjacency matrix of the bipartite realization $G$ where the columns are indexed by the vertices of $V$, numbered from left to right, and the rows are indexed by the vertices of $U,$ numbered from bottom to top. Hence the entry in row $i,$ column $j$ of the matrix will be written as $(j, i)$ and corresponds to the chord $(v_j , u_i).$ With some abuse of notation we also will use the word  ``chord" to refer to the matrix position as well. This
is nonstandard notation for the entries of a matrix, but matches the Cartesian coordinate system.  Then  let
\begin{displaymath}
\widehat M(X+Y-Z) = M_X + M_Y - M_Z.
\end{displaymath}
By definition each entry of an adjacency matrix is $0$  or $1$. Therefore only
$-1,0,1,2$ can be the entries of $\widehat M.$ An entry is $-1$ if the corresponding edge is
missing from both $X$ and $Y$ but it exists in $Z.$ The entry is $2$ if the corresponding edge is missing from $Z$ but exists in both $X$ and $Y.$ The entry is $1$ if the corresponding edge exists in all three graphs ($X,Y,Z$) or it is there only in one of $X$ and $Y$ but not in $Z.$ Finally it is $0$ if the corresponding edge is missing from all three graphs, or the edge exists in exactly one of $X$ and $Y$ and is also present in $Z.$ (Therefore if a chord denotes an existing edge in exactly one of $X$ and $Y$ then the entry corresponding to this chord is always $0$ or $1$.)
\begin{observ}\label{th:widehat}
Let $X, Y$ and $Z$ be some realizations of a bipartite degree sequence.
\begin{enumerate}[{\rm (i)}]
\item The row and column sums of $\widehat M(X+Y-Z)$ are the same as the row and column sums in $M_X$ (or $M_Y$ or $M_Z$).
\item If $Z$ is an intermediate realization between $X$ and $Y$ then $\widehat M(X+Y-Z)$ is another realization of the same degree sequence $($and all entries are $0$ or $1)$.
\end{enumerate}
\end{observ}

\bigskip\noindent Before we define some further notions we  introduce our main tool that  we will use later in this paper to illustrate different procedures in our current realizations.

Usually each cycle under processing is small comparing with the full graph, therefore we always consider a ``comfortably reordered'' adjacency matrix (in other words, we apply a suitable permutation on the vertices) such that the vertices forming the cycle  will be associated to an $\ell \times \ell$ submatrices of our adjacency matrices, and our figures will show only these submatrices. The positions $(1,1),\ldots, (\ell,\ell)$ form the {\bf main-diagonal} while the positions right above the main-diagonal as well as the rightmost bottom one (these are $(1,2), (2,3),\ldots,$ $(\ell-1,\ell)$ finally $(\ell,1)$)  form the {\bf small-diagonal}. (This placement was our goal using this numbering system for rows and columns. For example, the element $(1,2)$ corresponds to the chord $(v_1,u_2).$ If this is $1,$ then there is an edge there, otherwise the edge is missing.)

Now we introduce a new tool to give a slightly different view about this "central region".
This tool is the $\ell \times \ell$ matrix   $F_Z$: for a realization $Z$ where all chords are equal to the chords $G$ not completely within $V(C).$ In $V(C)$ (so at his central region) for $i,j = 1,\ldots,\ell$ we have:
$$
F_Z(j,i)=
\begin{cases}
M_Z(j,i) & \mbox{if } (j,i) \in \mbox{ main- or small-diagonals,} \\
\left [M_G+M_{G'} + M_Z\right ] (j,i) & \mbox{otherwise.}
\end{cases}
$$
 In that way in the main- and small-diagonal's elements are $0$ or $1$ while the others (the {\bf off-diagonal} entries) can be $0,1,2,3.$ There is an easy algorithm to construct $F_Z$ from the corresponding $\widehat M(G + G' - Z)$ and vice versa (please recognize that here we use $G$ and $G'$ instead of $X$ and $Y$): In the main-diagonal and in the  small-diagonal the zeros and ones must be interchanged. Outside of these diagonal entries $-1, 0, 1, 2$ of $\widehat M(G + G' - Z)$ become $1,0,3,2$ in $F_Z.$ (In case we need a second realization, similar to $Z$, we will denote it with $Z'.$)

Since $G$ and $G'$ coincide outside  the alternating cycle $C$ therefore the off-diagonal elements in $F_Z$ are odd when the edge exists in the actual $Z$ and even otherwise. When $Z=G$ then the main-diagonal entries are $1$ while the small-diagonal elements are $0.$ This matrix $F_Z$ will be used in our illustrating figures and also to conduct the construction of our canonical path system.

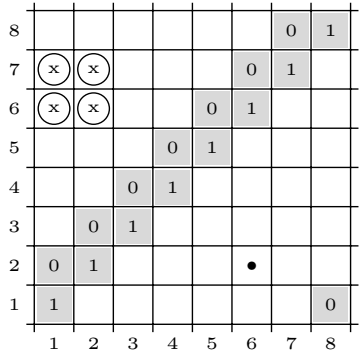
\begin{figure}[h]
\caption{A chord and its cousins}\label{fig1}
\center{\proba{
	\scalebox{1.3}{
	\begin{tikzpicture}[scale=0.4]
	\alap {8}
	\hh {6}{2}{$\bullet$}	\hhc {2}{6}{x}	\hhc {2}{7}{x}	\hhc {1}{6}{x}
	\hhc {1}{7}{x}
	\end{tikzpicture}
    }}}
\end{figure}

We are ready now to introduce the central notions of our proof:
\begin{definition}
The {\bf type} of a chord is $1$ if it is present in $G$, and $0$ otherwise. Note that a chord is present in $G$ if and only if it is present in $G'.$ Let $(v_\beta,u_\alpha)$ be a chord so $\delta\not\in \{\beta,\beta+1\}.$ A chord $(v_\beta,u_\alpha)$ is a {\bf cousin} of a chord  $(v_\delta, u_\epsilon)$,  if the other two corners of the submatrix, which is spanned by this position and  the chord are on the main- or on the small-diagonals of $F_Z$ (see Figure~\ref{fig1}). We can describe it with formulae as well: this chord $(v_\beta,u_\alpha)$ is a {\bf cousin} of a chord  $(v_\delta, u_\epsilon)$,  if $\alpha\not\in \{\beta, \beta+1\}$ and one of the following holds:
$$
\begin{cases}
\epsilon < \delta, & \alpha \in \{\delta,\delta+1\} \mbox{ and } \beta\in \{ \epsilon-1,\epsilon \}, \\
\epsilon >\delta, & \alpha \in \{\delta-1,\delta\} \mbox{ and } \beta\in \{ \epsilon,\epsilon+1 \}.
\end{cases}
$$
A chord $e$ is {\bf friendly} if at least one of its cousins has the same type as $e$ itself, otherwise it is {\bf unfriendly}. (Please recall that here ``chord" also refers to the position itself within the matrix therefore we also say that the position is friendly.)
\end{definition}
Now Figure~\ref{fig1} illustrates the cousins of the chord $(v_6,u_2)$ in the initial realization $Z=G.$ (They are $(v_1,u_6),$ $(v_1,u_7),$ $(v_2,u_6),$ finally $(v_2,u_7)$ and let's recall that the word chord indicates that the definition does not depend on the actual existence or non-existence of that edge.)

\noindent Before the next important definition we introduce a metric on pairs of positions of this matrix: $\|A,\bar A \|$ says how many steps are necessary to go from $A$ to $\bar A$ if
in every step we can move to a (horizontally or vertically) neighboring position, we cannot cross the main-diagonal, finally the position $(i,1)$  is neighboring to $(i,\ell)$ and analogously $(\ell, i)$ is neighboring to $(1,i)$.

\begin{definition}\label{def:friendly}
A sequence of pairwise distinct positions $A_1,\ldots,A_j$ is a {\bf friendly path} in $F_G$
if
\begin{enumerate}[{\rm (i)}]
\item each position is friendly (in the matrix $F_G$),
\item $\| A_h, A_{h+1}\|=1$,
\item the chords $e_1$ and $e_j$ - defined by the positions $A_1$ and $A_j$ - are shortest chords and the root of $e_1$ belongs to $G$ while the root of $e_j$ does not.
\end{enumerate}
\end{definition}

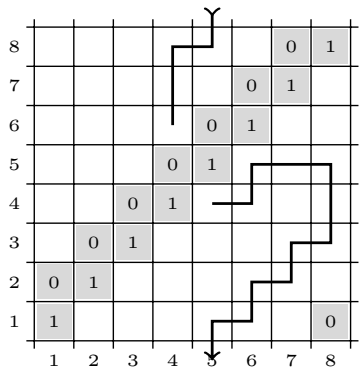
\begin{figure}[h]
\caption{A friendly path}\label{fig2}
\center{\proba{
	\scalebox{1.3}{
	\begin{tikzpicture}[scale=0.4]
	\alap {8}
\draw[thick,->](5.5,4.5)--(6.5,4.5)--(6.5,5.5)--(8.5,5.5)--(8.5,3.5)--(7.5,
3.5)-- (7.5,2.5)--(6.5,2.5)--(6.5,1.5)--(5.5,1.5)--(5.5,0.5);
    \draw [thick,>-](5.5,9.5)--(5.5,8.5)--(4.5,8.5)--(4.5,6.5);
	\end{tikzpicture}
    }}}
\end{figure}
\noindent A friendly path goes from the main-diagonal to the small-diagonal and it can be quite complicated, and it is important to remark that such a friendly path is NOT a path in a particular graph. Furthermore the friendly path is {\bf fixed} for the entire process determining the swap sequence from realization $G$ to realization $G'$, while the notions of chord or cousin apply for each matrix $F_Z$ along the swap sequence.

The name is justified by the image of the friendly path in the illustration of $F_G$, shown in Figure~\ref{fig2}. (It shows the path itself, but it does not show why the individual elements of the path are friendly.) The figures like this are not for illustration only: whenever we consider a friendly path we always work on the matrix itself.

\subsection{The case that a friendly path exists}\label{sec:friendly}
In this subsection we describe the construction of the path along this cycle in the case that a friendly path exists. Fix one friendly path: if there are more than one, then take, say, the lexicographically smallest one (relative to the subscripts of the positions). Let the chords of the existing friendly path correspond to the positions $A_1,\ldots,A_\Lambda$ where  $A_j= (a_j^1 ,a_j^2)$.

\medskip By definition our friendly path has the following properties: (i)
$a_j^1 \ne a_j^2$ and $a_j^1 + 1\ne a_j^2.$ (ii) $\|A_j,A_{j+1}\|=1,$ and finally
(iii) $A_1$ is at distance 1 from the main-diagonal, while $A_\Lambda$ is at
distance 1 from the small-diagonal.

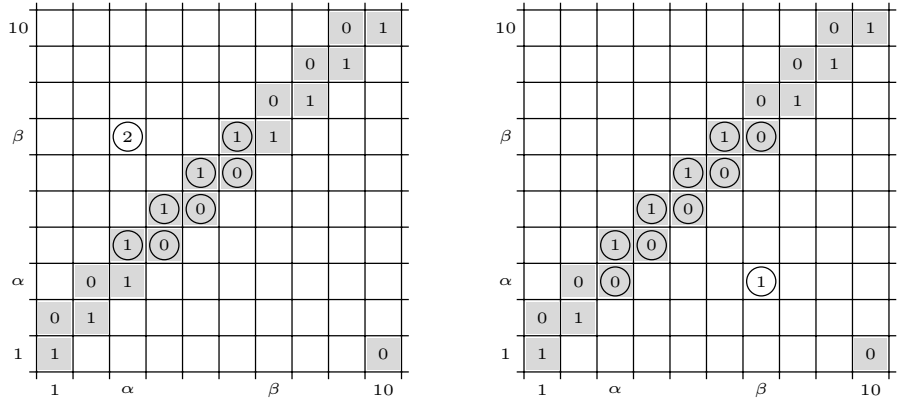
\begin{figure}[h]
\caption{An $(\alpha,\beta)$-OK matrix and an $(\beta,\alpha)$-KO
matrix}\label{fig3}
\center{\proba{\scalebox{1.2}{
\begin{tikzpicture}[scale=0.4]
    \hh{1}{0}{1}    \hh{3}{0}{$\alpha$}   \hh{7}{0}{$ \beta$} \hh{10}{0}{10}
    \hh{0}{1}{1} \hh{0}{3}{$\alpha$}  \hh{0}{7}{$\beta$} \hh{0}{10}{10}  \kalap {10}
    \hhc{3}{7}{2}  \rhhc341    \rhhc451    \rhhc561    \rhhc440    \rhhc550
\rhhc660    \rhhc671   \rhh{10}{1}{0}
\end{tikzpicture}
\qquad
\begin{tikzpicture}[scale=0.4]
    \hh{1}{0}{1}    \hh{3}{0}{$\alpha$}   \hh{7}{0}{$ \beta$} \hh{10}{0}{10}
    \hh{0}{1}{1} \hh{0}{3}{$\alpha$}  \hh{0}{7}{$\beta$} \hh{0}{10}{10}  \kalap
{10}
    \hhc{7}{3}{1}  \rhhc330    \rhhc341    \rhhc451    \rhhc561    \rhhc440
\rhhc550
    \rhhc660    \rhhc770    \rhhc671 \rhh{10}{1}{0}
\end{tikzpicture}
}}}
\end{figure}

\medskip\noindent First we introduce two new structures:
\begin{definition}
Let $1\le \alpha, \beta\le \ell$ with  $\beta\not\in \{\alpha,\alpha+1\}$. We say an $\ell\times
\ell$-matrix $F_Z=(m_{i,j})$ is {\bf $(\alpha,\beta)$-OK matrix} iff
\begin{enumerate}[{\rm (i)}]
\item $m_{\alpha,\beta}=2$,
\item $m_{i,i}=
\begin{cases}
0 & \mbox{for } i=\alpha+1, \alpha+2, \ldots ,\beta-1, \\
1 & \mbox{for } i=\beta, \beta+1, \ldots, \alpha-1, \alpha,
\end{cases}
$
\item $m_{i,i+1}=
\begin{cases}
1 & \mbox{for } i=\alpha, \alpha+1, \ldots,  \beta-2, \beta-1,\\
0 & \mbox{for } i=\beta, \beta+1, \ldots, \alpha-2, \alpha-1.
\end{cases}
$
\end{enumerate}
\end{definition}
\noindent (See the LHS of Figure~\ref{fig3}.) Please recall that the entry $2$ in $F_Z$ is
an edge which is missing from $Z$ but exists in both $G$ and $G'$ (the
off-diagonal entries are the same in $M_G$ and $M_{G'}$).

\begin{definition}
Let $1\le \alpha,\beta\le \ell$ with $\beta \not\in \{\alpha-1, \alpha\}$. We say an $\ell\times
\ell$-matrix $F_Z=(m_{i,j})$ is {\bf $(\beta,\alpha)$-KO matrix} iff
\begin{enumerate}[{\rm (i)}]
\item $m_{\beta,\alpha}=1$,
\item  $m_{i,i}=
\begin{cases}
0 & \mbox{for } i=\alpha, \alpha+1, \ldots, \beta-1, \beta, \\
1 & \mbox{for } i=\beta+1,\ldots, \alpha-1,
\end{cases}
$
\item $m_{i,i+1}=
\begin{cases}
1 & \mbox{for } i=\alpha, \alpha+1, \ldots, \beta-1, \\
0 & \mbox{for } i=\beta, \beta+1, \ldots, \alpha-2, \alpha-1.
\end{cases}
$
\end{enumerate}
\end{definition}
\noindent (See the RHS of Figure~\ref{fig3}.) Please recall that the entry $1$ in $F_Z$ is
an edge which exists in $Z$ but missing from both $G$ and $G'$.
\begin{lemma}\label{lm:OK}
Let $F_Z=(m_{i,j})$ be an $(\alpha,\beta)$-OK matrix and $m_{\alpha-1,\beta+2}=3$.
Assume that $F_{Z'}=(m'_{i,j})$ is an $(\alpha\!-\!1,\beta\!+\!2)$-OK matrix such
that
\begin{enumerate}[{\rm (1)}]
\item $m'_{\alpha,\beta}=3,$
\item $m'_{i,j}=m_{i,j}$ if $i\ne j$, $i+1\ne j$, and $(i,j)\ne (\alpha,\beta)$,
$ (\alpha-1,\beta+2)$.
\end{enumerate}
Then there exists an absolute constant  $\Theta$ such that one can transform $Z$ to $Z'$ by at most $\Theta$ swaps $($and, meanwhile,  transform $F_Z$ into $F_{Z'})$.
\end{lemma}
\proba{\scalebox{1.15}{\raisebox{2.0cm}{\hspace{-.6cm}$F_Z$}\!\!\!
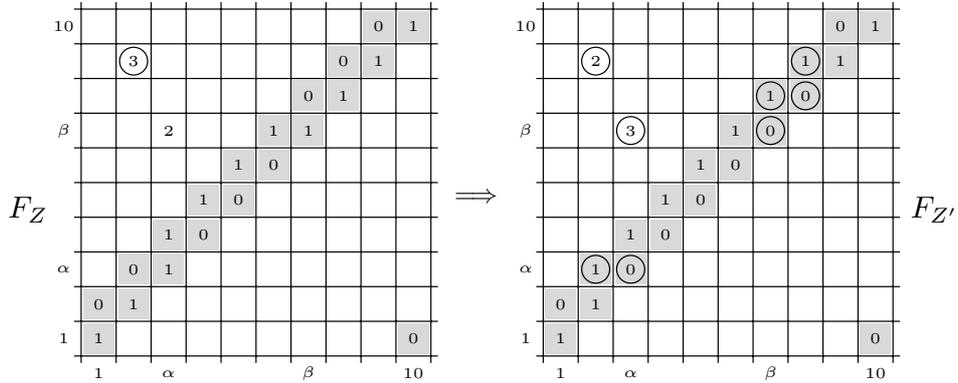
\begin{tikzpicture}[scale=0.4]
    \hh{1}{0}{1}    \hh{3}{0}{$\alpha$}   \hh{7}{0}{$ \beta$} \hh{10}{0}{10}
    \hh{0}{1}{1} \hh{0}{3}{$\alpha$}  \hh{0}{7}{$\beta$} \hh{0}{10}{10}  \kalap {10}
    \hh{3}{7}{2}    \rhh341    \rhh451    \rhh561    \rhh440    \rhh550
\rhh660    \rhh671    \hhc293 \rhh{10}{1}{0}
\end{tikzpicture}
    }
\raisebox{2.5cm}{$\Longrightarrow$}\!\!
\scalebox{1.15}{
\begin{tikzpicture}[scale=0.4]
    \hh{1}{0}{1}    \hh{3}{0}{$\alpha$}   \hh{7}{0}{$ \beta$} \hh{10}{0}{10}
    \hh{0}{1}{1} \hh{0}{3}{$\alpha$}  \hh{0}{7}{$\beta$} \hh{0}{10}{10}  \kalap {10}
    \hhc{3}{7}{3}  \rhh341    \rhh451    \rhh561    \rhh440    \rhh550
\rhh660    \rhh671
    \hhc292    \rhhc231    \rhhc891    \rhhc781    \rhhc770    \rhhc880
\rhhc330 \rhh{10}{1}{0}
\end{tikzpicture}
\raisebox{2cm}{$F_{Z'}$}
}}

\bigskip \Proof : It is enough to  observe that the symmetric difference of $Z$ and $Z'$ is a single alternating cycle. Indeed, in the next figure the entry $1$ indicates edges in $E(Z' - Z)$ and the entry $0$ indicates edges in $E(Z-Z').$ (The non-empty positions of this figure are the circled positions in the previous matrix $F_{Z'}$.)

\bigskip\noindent\hspace{-.5cm}
\begin{minipage}{7cm}
\center{$E(Z'-Z) \cup E(Z - Z')$}\\ \vspace{-.4cm}
\center{\proba{\scalebox{1.3}{
\begin{tikzpicture}[scale=0.4]
\alapalap {10}
    \hh23{1}    \hh33{0}    \hh290    \hh371    \hh770    \hh781    \hh880
\hh891
    \hh{3}{0}{$\alpha$}    \hh{0}{7}{$\beta$}
\end{tikzpicture}
}}}
\end{minipage}
\quad
\begin{minipage}{5cm}
Therefore\\
$\displaystyle(\alpha-1,\alpha), (\alpha,\alpha),  (\alpha,\beta),$\\
$\displaystyle(\beta,\beta), (\beta, \beta+1), (\beta+1, \beta+1),$\\
$\displaystyle(\beta+1, \beta+2), (\alpha-1,\beta+2)$ is \\
an alternating cycle of length 8.
\end{minipage}

\bigskip\noindent So the difference of the realizations lay in the subgraphs induced by $\bar V,$ which subset contains $8$ vertices. The subgraphs $Z[\bar V]$ and $Z'[\bar V]$ induced by $\bar V$ have the same (bipartite) degree sequence and they contain alternately the edges of the cycle. By Theorem~\ref{th:ryser} we know one of them can be transformed by swaps into the other one. Since the cycle contains  four-four vertices from both classes, and there are at most $12$ edges, therefore the canonical swap sequence (by  Corollary ~\ref{lm:kicsi})  is at most $2\times 12$  long therefore $\Theta = 24 $  is an upper bound  on the number of the necessary swaps. \qed

\medskip Clearly the same argument gives the following more general lemma.
\begin{lemma}\label{lm:uOK}
For each natural number $u$ there is a natural number $\Theta_u$ with the
following property: assume that $F_Z=\left (m_{i,j} \right )$ is an
$(\alpha,\beta)$-OK matrix and $m_{\alpha',\beta'}=3$ where
\begin{displaymath}
\big \|(\alpha,\beta);(\alpha',\beta')\big \| = u,
\end{displaymath}
furthermore $F_{Z'}=(m'_{i,j})$ is an $(\alpha',\beta')$-OK matrix such that
\begin{enumerate}[{\rm (1)}]
\item $m'_{\alpha,\beta}=3,$
\item $m'_{i,j}=m_{i,j}$ if $i\ne j$, $i+1\ne j$, and $(i,j)\ne (\alpha,\beta)$,
$ (\alpha',\beta')$.
\end{enumerate}
Then at most $\Theta_u$ swaps transform $Z$ into $Z'$  $($and along this $F_Z$ is transformed into $F_{Z'})$.
\end{lemma}
\Proof The only difference is that here the symmetric difference of $Z$ and $Z'$
is a cycle of length at most $2+2u$ which alternates between $Z$ and $Z'$. \qed

\medskip\goodbreak\noindent We also have the analogous general result for KO
matrices.
\begin{lemma}\label{lm:uKO}
For each natural number $u$ there is a natural number $\Theta'_u$ with the following property: assume that $F_Z=\left (m_{i,j} \right )$ is an $(\beta,\alpha)$-KO matrix and $m_{\beta',\alpha'}=0$ where
\begin{displaymath}
\big \| \left (\beta,\alpha);(\beta',\alpha'\right )\big \| = u,
\end{displaymath}
furthermore $F_{Z'}=\left (m'_{i,j}\right )$ is an $(\beta',\alpha')$-KO matrix such that
\begin{enumerate}[{\rm (1)}]
\item $m'_{\beta,\alpha}=0,$
\item $m'_{i,j}=m_{i,j}$ if $i\ne j$, $i+1\ne j$, and $(i,j)\ne (\beta,\alpha)$,
$ (\beta',\alpha')$.
\end{enumerate}
Then at most $\Theta'_u$ swaps transform $Z$ into $Z'$ $($and $F_Z$ is transformed into $F_{Z'})$.
\end{lemma}
\Proof The proof is very similar to the proof of Lemma~\ref{lm:uOK} which is left to the diligent reader. \qed

\begin{lemma}\label{lm:OK-KO}
Assume that $F_Z=(m_{i,j})$ is $(\alpha,\beta)$-OK matrix and $m_{\beta+2,
\alpha-1} =0$. Assume that $F_{Z'}=(m'_{i,j})$ is a $(\beta+2,\alpha-1)$-KO matrix
such that
\begin{enumerate}[{\rm (1)}]
\item $m'_{\alpha,\beta}=3,$
\item $m'_{i,j}=m_{i,j}$ if $i\ne j$, $i+1\ne j$, and $(i,j)\ne (\alpha,l)$, $
(\beta+2,\alpha-1)$.
\end{enumerate}
Then there exists a natural number $\Omega$ such that one can transform $Z$ into
$Z'$ by at most $\Omega$ swaps $($and $F_Z$ goes into $F_{Z'})$.
\end{lemma}
\proba{\scalebox{1.15}{\raisebox{2cm}{\hspace{-.6cm}$F_Z$}\!\!\!
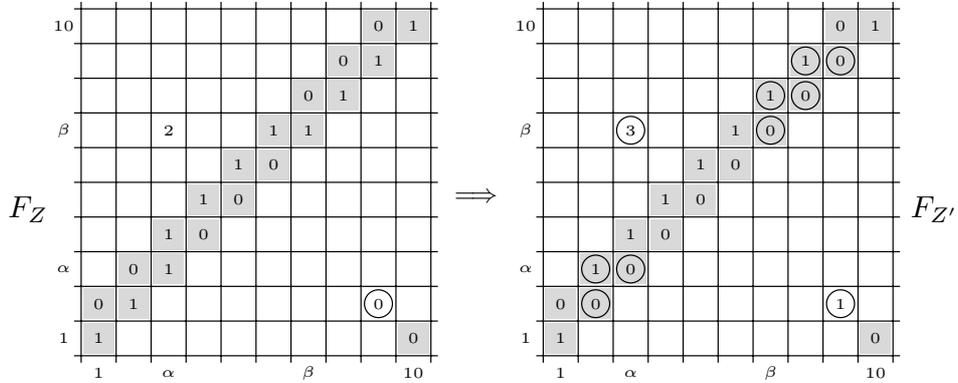
\begin{tikzpicture}[scale=0.4]
    \hh{1}{0}{1}    \hh{3}{0}{$\alpha$}   \hh{7}{0}{$ \beta$} \hh{10}{0}{10}
    \hh{0}{1}{1} \hh{0}{3}{$\alpha$}  \hh{0}{7}{$\beta$} \hh{0}{10}{10}  \kalap {10}
    \hh{3}{7}{2}    \rhh341  \rhh451    \rhh561    \rhh440    \rhh550    \rhh660
   \rhh671    \hhc920 \rhh{10}{1}{0}
\end{tikzpicture}
}
\raisebox{2.5cm}{$\Longrightarrow$}\!\!
\scalebox{1.15}{
\begin{tikzpicture}[scale=0.4]
    \hh{1}{0}{1}    \hh{3}{0}{$\alpha$}   \hh{7}{0}{$ \beta$} \hh{10}{0}{10}
    \hh{0}{1}{1} \hh{0}{3}{$\alpha$}  \hh{0}{7}{$\beta$} \hh{0}{10}{10}  \kalap {10}
    \hhc{3}{7}{3}    \hhc921
    \rhh341    \rhh451    \rhh561    \rhh440    \rhh550    \rhh660    \rhh671
    \rhhc220    \rhhc330    \rhhc231    \rhhc990    \rhhc891    \rhhc880
    \rhhc770    \rhhc781 \rhh{10}{1}{0}
\end{tikzpicture}
\raisebox{2cm}{$F_{Z'}$}
    }}
\Proof It is enough to  observe that the symmetric difference of $Z$ and $Z'$ is a single alternating cycle. Indeed, in the next figure values $1$ indicate edges in $E(Z'- Z)$ and values $0$ indicate edges in $E(Z- Z').$

\bigskip\noindent
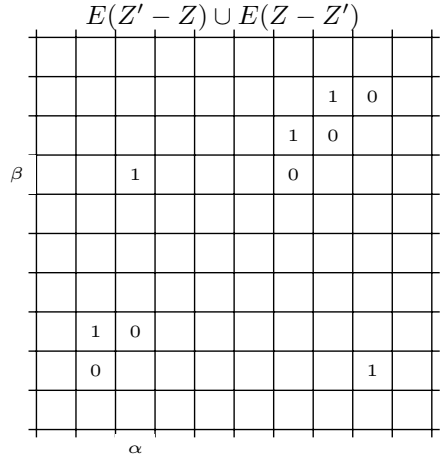
\begin{minipage}{7cm}
\center{$E(Z'-Z) \cup E(Z - Z')$}\\ \vspace{-.4cm}
\center{\proba{\scalebox{1.3}{
\begin{tikzpicture}[scale=0.4]
\alapalap {10}
    \hh220    \hh921    \hh990    \hh891    \hh880    \hh781    \hh770    \hh371
    \hh330    \hh231    \hh{3}{0}{$\alpha$}    \hh{0}{7}{$\beta$}
\end{tikzpicture}
}}}
\
\end{minipage}
\begin{minipage}{6cm}
Therefore\\
$\displaystyle(\alpha-1,\alpha-1),(\beta+2,\alpha-1), $\\
$\displaystyle(\beta+2,\beta+2), (\beta+1,\beta+2),$\\
$\displaystyle(\beta+1,\beta+1) (\beta,\beta+1), $\\
$\displaystyle(\beta,\beta),(\alpha,\beta),$\\
$\displaystyle(\alpha,\alpha), (\alpha-1,\alpha)$\\
is an alternating cycle of\\ length $10$.
\end{minipage}

\bigskip\noindent The proof goes like the proof of Lemma~\ref{lm:OK}: The
difference of the realizations lay in the subgraphs induced by  $\bar V,$ which subset contains $10$ vertices. The subgraphs $Z[\bar V]$ and $Z'[\bar V]$ induced by $\bar V$ have the same (bipartite) degree sequence and they contain alternately the edges of the cycle. By Theorem~\ref{th:ryser} we know one can be transformed by swaps into the other
one. Since the cycle contains  five vertices from both classes, and there are at
most $20$ edges, the number of the necessary swaps (by Corollary~\ref{lm:kicsi})
 is at most $2\times 20$   therefore there exists a constant upper bound $\Omega
\le 40 $ on the number of the necessary swaps. \qed

\begin{lemma}\label{lm:uOK-KO}
For each natural number $u$ there is a natural number $\Omega_u$ with the
following property: assume that $F_Z=(m_{i,j})$ is $(\alpha,\beta)$-OK and
$m_{\beta', \alpha'}=0$ where
\begin{displaymath}
\big \|(\alpha,\beta);(\alpha',\beta')\big \| = u,
\end{displaymath}
and $F_{Z'}=(m'_{i,j})$ is a $(\beta',\alpha')$-KO matrix such that
\begin{enumerate}[{\rm (1)}]
\item $m'_{\alpha,\beta}=3,$
\item $m'_{i,j}=m_{i,j}$ if $i\ne j$, $i+1\ne j$, and $(i,j)\ne (\alpha,\ell)$,
$ (\beta',\alpha')$.
\end{enumerate}
Then at most $\Omega_u$ swaps transform $Z$ into $Z'$ $($and  $F_Z$ into $F_{Z'})$.
\end{lemma}
\Proof Similar to Lemma~\ref{lm:uOK}. \qed

\bigskip\noindent  Now using our friendly path we are going to define a sequence
of OK- and KO-matrices, such that we can achieve the required edge changes in
$G$ obtaining $G'$ along this sequence, using operations described in the previous
Lemmas. At first we define a new sequence $A_1',\ldots,A_\Lambda'$ from
$A_1,\ldots,A_\Lambda$ in the following way:
\begin{equation}\label{eq:vesszo}
A_i'=
\left\{
  \begin{array}{ll}
    A_i, & \hbox{if } F_G(A_i)=0, \\
    \mathrm{Cousin}(A_i), & \hbox{if } F_G(A_i)=3,
  \end{array}
\right.
\end{equation}
where $\mathrm{Cousin}(A)$ denotes one of the  cousins of $A$. If there are more than one positions of the same type among the corresponding positions, then we choose the lexicographically-least one. We will use the following notation: the {\em mirror image} of the position $(\alpha,\beta)$ to the main-diagonal is $\mathrm{Mirror} (\alpha,\beta) =(\beta,\alpha)$.
\begin{observ}\label{lm:dist}
By definitions,
\begin{enumerate}[{\rm (i)}]
\item if $F_G(A_i)=F_G(A_{i+1})$ then $\left \|A'_i, A'_{i+1} \right \| \le 3,$
\item if $F_G(A_i)\ne F_G(A_{i+1})$ then $\left \|\mathrm{Mirror}(A'_i),
A'_{i+1}\right \| \le 3.$
\end{enumerate}
\end{observ}
\begin{definition}
We define the matrix sequence $F_G=L_0,L_1,\ldots,L_\Lambda,L_{\Lambda+1}=F_{G'}$ and the corresponding realizations $Z_1,\ldots, Z_\Lambda$, where $L_i = F_{Z_i}$ for each $i$ as follows:

The matrix $L_i$ ($i=1,\ldots,\Lambda$) is defined from the matrix $L_{i-1}$ by the formulae:
\begin{displaymath}\label{eq:matrix-sec}
L_i =
\left\{
  \begin{array}{ll}
    \hbox{the }(A'_i)\hbox{-OK matrix,} & \hbox{if } L_{i-1} (A_i)=3, \\
    \hbox{the }(A'_i)\hbox{-KO matrix}, & \hbox{if } L_{i-1} (A_i)=0.
  \end{array}
\right.
\end{displaymath}
Here all positions $(u,v)$ which are NOT determined by the definitions of the
OK- and KO-matrices satisfy $L_i(u,v)=L_{0}(u,v).$ \qed
\end{definition}
\noindent It is quite clear that $(\Lambda -1)$ consecutive applications of (the
appropriate) Lemmas~\ref{lm:OK} - \ref{lm:uOK-KO} will take care the definition
of the required swap sub-sequences between $L_1$ and $L_\Lambda.$ However, the
swap-sequence transforming $L_0$ into $L_1$ furthermore the one transforming
$L_\Lambda$ into $L_{\Lambda+1}$ require special considerations:
\begin{itemize}
\item If $L_0(A_1)=3$ then there are two possibilities - depending on the
position of the $\mathrm{Cousin}(A_1)$. (The squares denoted with dashed lines
contain the possible positions of friendly cousins.)\\

Case I:\\
\medskip\scalebox{1.6}{\proba
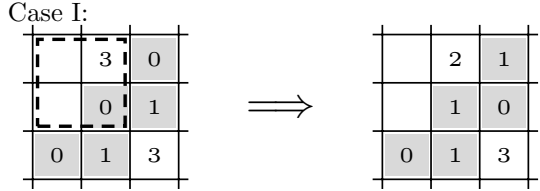
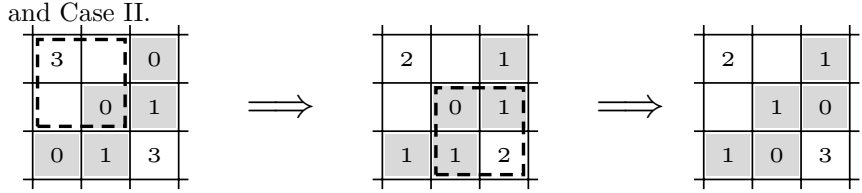
{
\begin{tikzpicture}[scale=0.4]
	\alapalap {3}
    \rhh110    \rhh220    \rhh330    \rhh211    \rhh321    \hh313    \hh233
    \draw[thick,\haha]  (1.1,2.1) rectangle (2.9,3.9);
\end{tikzpicture}
\quad \raisebox{.6cm}{$\Longrightarrow$} \quad
\begin{tikzpicture}[scale=0.4]
	\alapalap {3}
    \rhh110    \rhh221    \rhh331    \rhh211    \rhh320    \hh313    \hh232
\end{tikzpicture}
}}\\
and Case II. \\
\scalebox{1.6}{\proba{
\begin{tikzpicture}[scale=0.4]
	\alapalap {3}
    \rhh110    \rhh220    \rhh330    \rhh211    \rhh321    \hh313    \hh133
    \draw[thick,\haha]  (1.1,2.1) rectangle (2.9,3.9);
\end{tikzpicture}
\quad \raisebox{.6cm}{$\Longrightarrow$} \quad
\begin{tikzpicture}[scale=0.4]
	\alapalap {3}
    \rhh111    \rhh220    \rhh331    \rhh211    \rhh321    \hh312    \hh132
    \draw[thick,\haha]  (2.1,1.1) rectangle (3.9,2.9);
\end{tikzpicture}
\quad \raisebox{.6cm}{$\Longrightarrow$} \
\begin{tikzpicture}[scale=0.4]
	\alapalap {3}
    \rhh111    \rhh221    \rhh331    \rhh210    \rhh320    \hh313    \hh132
\end{tikzpicture}
}}
\item If, however, $L_0(A_1)=0$ then there is only one case:\\
\scalebox{1.6}{\proba
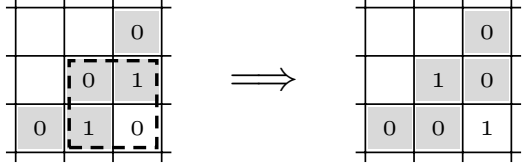
{
\begin{tikzpicture}[scale=0.4]
	\alapalap {3}
    \rhh110    \rhh220    \rhh330    \rhh211    \rhh321    \hh310
    \draw[thick,\haha]  (2.1,1.1) rectangle (3.9,2.9);
\end{tikzpicture}
\quad \raisebox{.6cm}{$\Longrightarrow$} \quad
\begin{tikzpicture}[scale=0.4]
	\alapalap {3}
    \rhh110    \rhh221    \rhh330    \rhh210    \rhh320    \hh311
\end{tikzpicture}
}}
\end{itemize}
The connecting swap-sequence from the matrix $L_\Lambda$ to $L_{\Lambda+1}$
(which is $F_{G'}$) can be defined analogously to the previous one. This completes the definition of the canonical path $\Gamma(X,Y,s).$

Next we will analyze the behavior of the current matrices $\widehat M(G,G',Z)$
along these sub-sequences. At first we consider those $Z$'s which correspond to
matrices $L_i.$

Let $M$ be an integer matrix and let $M'$ be a $2\times 2$ submatrix of it. If
we add $1$'s to the values of the positions of one diagonal in $M'$ and $-1$'s
to the values of the positions of the other diagonal, then the acquired matrix
has the same row and column sums as $M$ had. Such  an operation is called a {\bf
switch}. When our matrix $M$ is the adjacency matrix of a degree sequence
realization, then any swap clearly corresponds to a switch of that matrix. We say that the
two matrices are in {\bf switch-distance} $1$ from each other. It is clear that bounded switch-distance between two matrices also means bounded Hamming distance  between them (as it was required in (F)(d)).

The following lemma is an auxiliary result, which help us to handle the numbers of different paths (in our canonical path system) which cover the same edge. It has no role in the definition of our path system, but it helps to show that this path system obeys the rules outlined in (A) -- (F) is Section \ref{sec:general_new}.
\begin{lemma}\label{lm:switch}
For $i = 1, . . . , \Lambda$ there exist  realizations $G_1,\ldots, G_\Lambda$ in $V(\G)$ for which $M_{G_i}$ is in switch-distance $1$ from the matrix $\widehat M(G+G'-Z_i)$ for  $i = 1, . . . , \Lambda$.
\end{lemma}
\Proof We show here the statement for such an $L_i$  where $L_i(A_i)=0$ therefore
$L_i$ itself is an $(A'_i)$-KO matrix, and where - by definition - $A_i=A'_i$ (the other case is similar). Due to the definitions $A_i$ originally is not an edge either in $G$ or in $G'.$ It belongs to the friendly path, therefore we also know that  $F_G (\mathrm{Cousin} (A_i))=F_{G'}(\mathrm{Cousin}(A_i))=0$ hold. In $L_i$ this value is $1$, so $A_i$ is an edge in $Z_i.$ Therefore $L_i$ which is $=F_{Z_i}$ looks like the matrix to the left in the following figure (the circled element is the cousin of $A_i$). The corresponding $\widehat M(G+G'-Z_i)$ is shown on the right hand side:

\medskip\hspace{2.5cm}$F_{Z_i}$ \hspace{4.65cm}$\widehat M(G+G'-Z_i)$

\raisebox{2cm}{\ \vspace{3cm}} \scalebox{1.6}{
\begin{tikzpicture}[scale=0.4]
    \alap {6} \hh 521 \hhc 150  \rhh220 \rhh 550
\end{tikzpicture}
\qquad
\begin{tikzpicture}[scale=0.4]
    \alap {6} \hh 52{-1} \hhc 150
    \rhh 110 \rhh 660   \rhh 330 \rhh 440 \rhh 121  \rhh 231 \rhh 341 \rhh 451 \rhh 561 \rhh 611 \draw (5,2)--(5,3);  \draw (6,2)--(6,3);
\end{tikzpicture}
}
\\
It is clear that adding $1$ to the values of the positions $A_i$ and $\mathrm{Cousin} (A_i)$ of  $\widehat M(G+G'-Z_i)$  and subtracting $1$ from the other two corners of the spanned submatrix constitutes the required switch. \qed \\
( In the figures above $A_i$ is $(5,2)$.  Here one can also recall that outside our $\ell \times \ell$ submatrix every entry  is $0$ or $1$ and after the switch the same applies inside the submatrix. Therefore, due to the row- and column-sum conditions, the acquired matrix is a realization indeed.)
\begin{lemma}\label{lm:friendly-main}
The realization $G$ can be transformed into the realization $G'$ through
realizations $Z_i$ ($i=1,\ldots,\Lambda$) in such a way that the lengths of the swap sub-sequences leading from each $Z_i$ to $Z_{i+1}$ $($where $0=1,\ldots,\Lambda)$ can be bounded from above by the absolute constant $\max\{\Theta_3, \Theta'_3,$ $\Omega_3\}.$ In this process, each arisen matrix $\widehat M(G+G'-Z_i)$ is within a constant  switch-distance from some vertex in $V(\G)$ (that is some realization of the bipartite degree sequence).
\end{lemma}
\Proof By Observation~\ref{lm:dist} for each $i$ the positions $A_i'$ and $A'_{i+1}$ or $\mathrm{Mirror}(A_i')$ and $A'_{i+1}$  are at most distance $3.$ Therefore for each $i$ (where $i=2,\ldots,\Lambda$) the corresponding process chosen among Lemma \ref{lm:uOK}, Lemma \ref{lm:uKO} and Lemma \ref{lm:uOK-KO} will describe the desired swap sub-sequences. The length of any such swap-subsequence is bounded from above by $\max\{\Theta_3, \Theta'_3, \Omega_3\}.$

Furthermore when in the process the  current realization $Z_i$ corresponds to an
$F_{Z_i}=L_i,$ then Lemma~\ref{lm:switch} applies, and  matrix $\widehat M(G+G'-Z_i)$ has switch-distance $1$ from the adjacency matrix of some realization $\in V(\G)$.

Let now $Z$ be a realization in the process, say, on the path between the matrices $L_i$ and $L_{i+1}$: then $\widehat M(G+G'-Z_i)$ can be transformed through swaps into $\widehat M (G+G'-Z_{i+1})$ (assume, this end is the closer one to $Z$). As we know all swaps are specialized switches, and they keep the row and column sums. Combining this with the previous paragraph, we have for every $Z$ that $\widehat M(G+G'-Z)$ is at most $\left \lceil \frac{1}{2} \max\{\Theta_3, \Theta'_3, \Omega_3\}\right \rceil +1$ switch distance from some realization $\in V(\G).$ \qed

\bigskip\goodbreak
\centerline{\bf \large Key problem}

\medskip\noindent\fbox{\parbox{12cm}{
One can say that we are very close to proving the rapidly mixing property of our
Markov process on all bipartite degree sequences: we should prove, that in the
case when there exists a friendly path from $G$ to $G'$ then for each
intermediate $Z$ the matrix $\widehat M(X+Y-Z)$ is in a constant distance from some
realization $\in V(\G).$ If we can manage this then we must handle  the cases
when there are no friendly paths. It is somewhat surprising that this second
requirement can be satisfied successfully (as it will be shown in Subsection
\ref{sec:nofriendly}).

\quad However, we cannot manage to prove the first requirement. The problem is
the following: we can try to repeat the proof of Lemma~\ref{lm:friendly-main},
but, unfortunately, it is not true anymore that for each graph $Z$,
corresponding to a particular matrix $L_i$, the matrix $\widehat M(X+Y-Z)$ is also
in distance $1$ from some realization in $V(\G).$

\quad In the realizations $G$ and $G'$ all chords have the same types, but this
is not the case for realizations $X$ and $Y.$ The edges in $E(X-Y) \cup E(Y-X)$
belong to only one of them. Therefore if a swap turns an entry to $2$ in $\widehat
M(G+G'-Z)$ then this entry originally was $1$: the edge belonged to $G$ and $G'$
and $Z$ as well. Therefore its cousin bears the entry $1$ (also belonged to $G$
and $G'$ and $Z$ as well). So this entry was appropriate to perform a switch to
turn the matrix under investigation into the adjacency matrix of a realization.
However, if the cousin entry is $0$ in $\widehat M(X+Y-Z)$ (this edge belongs only
to one of realizations $X$ and $Y$, say, it belongs to $X$ only), then the
required switch cannot be performed. (The value $-1$ can cause a similar
problem and can be handled similarly as this case.)

\quad A good solution for this particular problem would probably end up in a
complete proof of the rapidly mixing property.
    }}

\bigskip\noindent The following observation is enough to handle the
switch-distance problem for $\widehat M(X+Y-Z)$  in  half-regular bipartite degree
sequences. Recall, a bipartite degree sequence $(\mathbf{a}, \mathbf{b})$ is {\em half-regular} if in $\mathbf{a}$ all degrees are the same, while the entries in $\mathbf{b}$ can be anything.
\begin{lemma}\label{lm:semi-reg}
Assume that our bipartite degree sequence  $(\mathbf{a}, \mathbf{b})$ is
half-regular and  the matrix $F_G$ under investigation contains a friendly path. Then the statement of Lemma~{\rm \ref{lm:friendly-main}} applies for the matrices $\widehat M(X+Y-Z_i)$ as well.
\end{lemma}
\Proof We follow the proof of Lemma~\ref{lm:friendly-main}. To do so the only
requirement is to show (somewhat loosely) that
the matrices $\widehat M(X+Y-L_i)$ are in a constant switch-distance from the
adjacency matrix of some realizations. As we know any of these matrices contains
exactly one entry of value different from $1$ and $0$. So consider a particular
$L_i$ and assume that this ``extra'' value in this case is a $2$. If the switch,
described in the proof of Lemma~\ref{lm:switch}, is also a possible switch in
$\widehat M(X+Y-Z)$ then we are done. If this not the case then the entry (with
value $1$ in matrix $\widehat M(G+G'-L_i)$) has value $0$ in $\widehat M(X+Y-L_i)$. (In
this case, as we discussed it previously, the corresponding edge is missing from
$Y.$) Let this corresponding edge be $(u,v)$, then this entry in $\widehat
M(X+Y-L_i)$ is $0.$ Since the column sums are fixed in these matrices, they are
the same (and equal to entries in $\mathbf{a}$).

Now vertex $v$ has degree at least $2$ (it is a vertex on cycle $C$ and it also end point of at least one chord of $C$ in $X$). Therefore the column $v$ contains some $1$s. One of them is $(w,v)$ (this $w$ cannot be the row of the $2$, since the entry there is $0$ due that it belongs to the originally intended switch). Now by the pigeonhole principle (since all row sums are the same) there is a column $z$ such that $\widehat M(w,z)=0$ and $\widehat M(u,z)=1$. Therefore the $u,w;v,z$ switch (actually this is a swap) will change $\widehat M(u,v)$ into $1$, and now the original switch finishes the job. The matrix $\widehat  M(X+Y-L_i)$ is in switch-distance at most 2 from the adjacency matrix of some realization. \qed

\subsection{The case that no friendly path exists} \label{sec:nofriendly}
In the previous subsection we discussed the situation when -- processing one by one the cycles in the canonical decomposition of the symmetric difference -- the cycle under investigation possesses a friendly path. All definitions, statements, reasonings were valid for any arbitrary bipartite degree sequence -- except the situation described in the  Key Problem and in Lemma \ref{lm:semi-reg} where we have to use the half-regularity condition.

Here we discuss the case where there exists no friendly path in the cycle under investigation. Nothing that we define here, state here or prove here requires the half-regularity condition.
So here our general assumptions are: we have realizations $G$ and $G'$ of the same (arbitrary) bipartite degree sequence, where the symmetric difference of the two edge set forms exactly one cycle, which, in turn does not possesses a friendly path.

Our plan is this: at first we show that the non-existence of the friendly paths
yields a strong structural property of the matrix $F_G.$ Using this property we
can divide our problem into two smaller ones, where one of the smaller matrices
possesses a suitable friendly path. So we can solve our original problem in a
recursive manner.

This recursive approach must be carried out with caution: a  careless ``greedy"
algorithm can increase the switch-distances very fast. We will deal with this
problem using a simple ``fine tuning" (which is described at the end of this
subsection).

We start with some further notions and notations.
\begin{definition}
In an $\ell \times \ell$ matrix the {\em sequence} of positions $(i+1, i-1), (i+2,i-2),\ldots, (i+\lfloor \ell /2 \rfloor -1,i- \lfloor \ell /2 \rfloor +1)=(j-1,j+1)$ form the  $\mathit{i}$th {\bf   down-line} of the matrix. (The arithmetic operations are thought to be considered modulo $\ell$, that is, for example,  $1-3=\ell - 2$. Therefore if the down-line reach the edge of the matrix at position, say, $(\ell, k)$ then the next position is $(1,k-1).$ Similarly, if the position on the edge is $(k,1)$ then the next position is $(k+1,\ell).$ If $2i>\ell$ then the first case applies, in  case of $2i < \ell$ the second case applies. Finally if, by chance, $2i=\ell$ then the positions in questions are $(\ell,1)$ and $(1,\ell).$  Analogously the {\em sequence} of positions  $(i-1, i+1), (i-2,i+2),\ldots,(i-\lceil \ell /2 \rceil + 1,i+\lceil \ell /2 \rceil -1)$ form the  $\mathit{i}$th {\bf up-line} of the matrix. (Let us mention that in case of even $\ell$ the length of the up-lines and the down-lines are equal. However, in case of odd $\ell$ the down-lines are longer with one position.)
\end{definition}
Since the lines are sequences therefore by definition they have orientations along which the algorithm will traverse them. Also by definitions the previous $\mathit{i}$th down-line and the $\mathit{j}$th up-line (for some $j$) in case of even $\ell$, as sets, are equal. However, as sequences, they are of course different.
\begin{definition} A set $T$ of positions of an   $\ell \times \ell$ matrix is called {\bf rook-connected} if a chess rook, staying inside $T$, can visit all elements of $T$. Here the chess rook is allowed to wrap around cyclically on the rows and columns (that is the rook is moving on a torus). We use the expression {\bf king-connected} analogously.
\end{definition}

The following lemma is a well-known version of the classical Steinhaus lemma
(see \cite{stein}).
\begin{lemma}\label{lm:stein}
Assume that the off-diagonal positions of an  $\ell \times \ell$ matrix are arbitrarily colored white and black. Then either the rook  has a white path which starts at distance 1 from the small-diagonal and ends at distance 1 from the main-diagonal, and avoids both diagonals, or there is a king-connected set $T$ of black positions  which intersects all rook's paths from the main-diagonal to the small-diagonal.  \qed
\end{lemma}
We use the previous result without proof. The set $T$, which was identified in
the previous lemma, will be called a {\em Steinhaus set}.
\begin{definition}
The {\bf cousin-set} $\mathfrak{C}(u,v)$ is the set of the off diagonal cousins
of the position $(u,v)$.  If $T$ is a set of positions, then the cousin set
$\mathfrak {C} (T)$ is defined as $\bigcup\{\mathfrak{C}(e):e\in T\}$.
\end{definition}
\begin{lemma}\label{lm:egyszinu}
Assume that in the matrix $F_G$ there is a king-connected set $T$ of unfriendly
positions. Then the  cousin-set $\mathfrak{C}(T)$ is rook-connected and type of all its positions are the same. All positions in $T$ have the opposite type.
\end{lemma}
\Proof W.l.o.g. we may assume that a position $P$ in $T$ has type $0,$ then all positions in its cousin-set must have type $1$. However, for each other position $P'$ in $T$, which can be reached from $P$ in one king step, the cousin sets $\mathfrak {C}(P)$ and  $\mathfrak{C}(P')$ have at least one common position. Therefore the neighboring cousin sets are rook-connected, furthermore all types in those two cousin sets must be the same (1), therefore both positions $P$ and $P'$ have the same type ($0$) as well. \qed
\begin{lemma}\label{lm:intersect}
Let $T$ be a Steinhaus set in $F_G$, then its cousin-set  $\mathfrak{C}(T)$  intersects all down-lines and up-lines.
\end{lemma}
\Proof Actually we can prove more: namely that any king-path from the main-diagonal to the small-diagonal intersects the cousin-set $\mathfrak{C}(T)$. Now the statement is equivalent with Lemma \ref{lm:stein} if we rotate the chess-board with 90 degree. Finally it is clear, that every down- and up-line forms a required king-path. \qed
\begin{lemma}\label{lm:same-state}
We assume that in the matrix $F_G$ there is no friendly path. Then for each $i$ $(i=1,\ldots, \ell)$ there exists a
$\mathbf{t} \in \{0,1\}$ and a pair of indices $j,j' \in \{1,\ldots \lceil \ell /2 \rceil -1  \}$ such that one of the following holds:
\begin{itemize}
\item entries  $(i+1,i-1),\ldots, (i+j, i-j)$ and $(i-j',i+j')$ have   the same type $\mathbf{t}$, furthermore  the entries $(i-1,i+1),   \ldots, (i-j'+1, i+j'-1)$ have the type $1-\mathbf{t}$, and all   entries belong to a down- or up-line;
\item entries  $(i-1,i+1),\ldots, (i-j, i+j)$ and $(i+j',i-j')$ have   the same type $\mathbf{t}$, furthermore  the entries $(i+1,i-1),   \ldots, (i+j'-1, i-j'+1)$ have the type $1-\mathbf{t}$, and all   entries belong to a down- or up-line.
\end{itemize}
\end{lemma}
\Proof Assume for a contradiction that there is no such $j$ for a particular $i.$ W.L.O.G. we
may assume, that $F_G(i+1,i-1)=0.$ Then, by the assumption, $F_G(i-1,i+1)= F_G(i-2,i+2)=1$ must hold. Then, again by our assumption, $F_G(i+2,i-2)=0$ must hold, etc. All entries along the down-line are $0$, while all entries along the up-line must be $1.$ However both lines intersect (see Lemma~\ref{lm:intersect}: both lines can be traversed by a chess king) the cousin-set $\mathfrak{C}(T)$ of the Steinhaus set $T$.  But, by Lemma~\ref{lm:egyszinu}, all its entries have the same type. A contradiction. \qed
\begin{corollary}\label{lm:van-friendly}
If conditions of Lemma~\ref{lm:same-state} hold, and $j'\ge 2$ for a particular
$i$, the submatrix spanned by $(i+j,i-j)$ and $(i-j,i+j)$ contains at least one friendly path. (Let us recall that the bottom-right position $(i+j,i-j)$ belongs to the small-diagonal by definition.)
\end{corollary}
\Proof We argue by contradiction: assume that the submatrix does not contain a friendly path. Then - due to Lemma~\ref{lm:stein} - it contains a Steinhaus set $T$. Due to  Lemma ~\ref{lm:egyszinu}, in its cousin-set $\mathfrak{C}(T)$ - which intersects all down- and up-lines - all positions have the same type. But this contradicts to the fact, that in the $i$th down-line all positions have type $\mathbf{t}$, while in the $i$th up-line all positions have type $1-\mathbf{t}$. A contradiction, again. \qed

\medskip\noindent That finishes the preliminaries that are needed to describe our recursive algorithm,  which is essentially a divide and conquer approach. Due to the previous fact here we should handle separately two possibilities: when $j' =1 $ (and $j=1$ as well) and when $j' \ge 2.$ For sake of simplicity we will assume that in our cycle the first condition described in Lemma \ref{lm:same-state} holds. We start with the

\bigskip \noindent {\bf First possibility}: assume that, for a particular $i,$ we have $j'=1.$ Then we also have $j=1.$ We should take care of two cases:

\medskip\noindent{\bf Case 1:}  If both $F_G(i+1,i-1)$ and $F_G(i-1,i+1)=3$ (that is both chords belong to both $G$ and $G'$) then we are in an easily handleable situation: at first we swap the quartet $u_i, u_{i+1}; v_i,v_{i+1}$.  (The dashed square in our illustration. Here we use the matrix $F_G$.) The entries $(i-1,i), (i,i)$ and $(i,i+1)$ have  the required types. However, entry $F'_G(i-1,i+1)=2$ therefore during the procedure  we should take care to change it back to its original value.

\medskip
\noindent \scalebox{1.5}{\proba
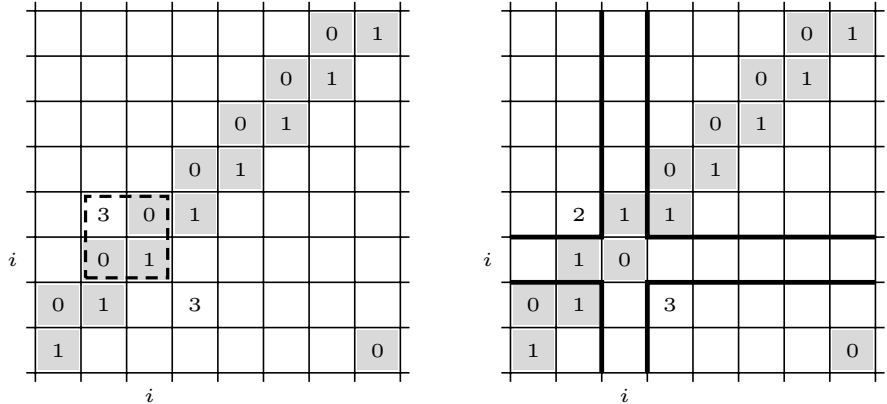
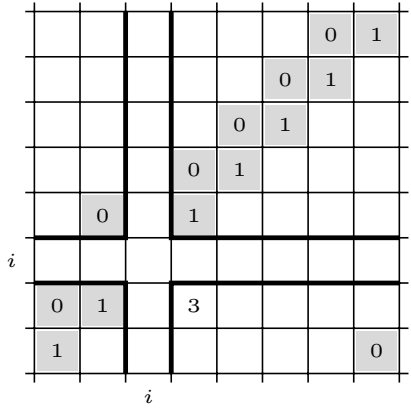
{
\begin{tikzpicture}[scale=0.4]
   \hh 30{$\scriptsize i$}  \hh 03{$\scriptsize i$} \kalap {8}
    \hh 423 
    \hh 243 
    \rhh{8}{1}{0}
    \draw[thick,\haha]  (2.1,3.1) rectangle (3.9,4.9);
\end{tikzpicture}
\quad
\begin{tikzpicture}[scale=0.4]
    \hh 30{$\scriptsize i$}  \hh 03{$\scriptsize i$} \kalap {8}
    \hh423
    \hh242
    \rhh330
    \rhh231
    \rhh341\rhh{8}{1}{0}
    \draw[very thick]  (1,3)--(3,3)--(3,1);
    \draw[very thick]  (1,4)--(3,4)--(3,9);
    \draw[very thick]  (9,4)--(4,4)--(4,9);
    \draw[very thick]  (4,1)--(4,3)--(9,3);
\end{tikzpicture}
}}

\noindent
\begin{minipage}{6cm}
\noindent \scalebox{1.5}{
\begin{tikzpicture}[scale=0.4]
    \hh 30{$\scriptsize i$}  \hh 03{$\scriptsize i$} \kalap {8}
    \hh423
    \rhh240
    \hh 33{\phantom{0}}  
    \hh 23{\phantom{1}}   
    \hh 34{\phantom{1}}   
    \rhh{8}{1}{0}
    \draw[very thick]  (1,3)--(3,3)--(3,1);
    \draw[very thick]  (1,4)--(3,4)--(3,9);
    \draw[very thick]  (9,4)--(4,4)--(4,9);
    \draw[very thick]  (4,1)--(4,3)--(9,3);
\end{tikzpicture}
}
\end{minipage}
\quad
\begin{minipage}{5.2cm}
The remaining subproblem, indicated by thick black lines,  fortunately is already in the required form. Indeed: its main-diagonal contains only $1$s, while its small-diagonal is full
with $0$s. (We have to keep it in our mind that the shown matrix of the remaining smaller subproblem is $F_Z$ where the element of the main- and small-diagonals came from $Z.$)
\end{minipage}

\noindent Denote the alternating cycle of this smaller problem by $C'$. The (recursive) subproblem $C'$ may contain a friendly path which will process it completely in one step, and will  switch the value of $F_G(i-1,i+1)$ automatically back to $1$.  If, however, it does not contain a friendly path, then the recursive procedure can use any down- and up-lines, including  $(i-1,i+1)$ (see Corollary~\ref{lm:intersect}), therefore we can take care that this switch-back will happen in the next recursion.

It is important to recognize that matrices $\widehat M(G+G'-Z)$ and $\widehat M(X+Y-Z)$
may contain $2$ at the position  $(i-1,i+1)$. Fortunately this ``problematic'' entry will be present only along one recursive step. Furthermore this entry will increase the switch-distance of the current $\widehat M$ by at most one:  the positions $(i-1,i), (i,i)$ and $(i,i+1)$ (outside of our subproblem), provide a suitable switch to handle the entry $2$ at position $(i-1, i+1)$ .

\bigskip\noindent{\bf Case 2:} Now we have $F_G(i+1,i-1)=0$ and $F_G(i-1,i+1)=0$. Here
we perform two swaps, the places of the swaps are denoted (shown below) with dashed squares:

\noindent
{\center\scalebox{1.3}{\proba
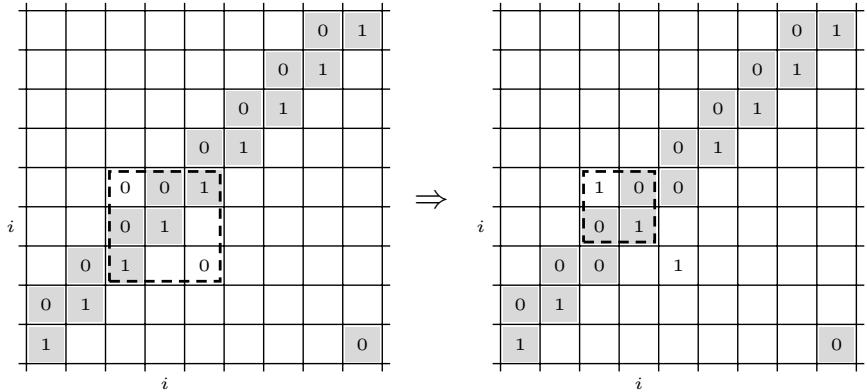
{
\begin{tikzpicture}[scale=0.4]
    \hh 40{$\scriptsize i$} \hh 04{ $\scriptsize i$} \kalap {9}
    \hh530    \hh350
    \draw[thick,\haha]  (3.1,3.1) rectangle (5.9,5.9);\rhh{9}{1}{0}
\end{tikzpicture}
\ \raisebox{2cm}{$\Rightarrow$}
\begin{tikzpicture}[scale=0.4]
    \hh 40{$\scriptsize  i$}  \hh 04{$\scriptsize i$}  \kalap {9}
    \hh531    \hh351    \rhh330    \rhh550
    \draw[thick,\haha]  (3.1,4.1) rectangle (4.9,5.9);\rhh{9}{1}{0}
\end{tikzpicture}
	}}}

\noindent  If $\widehat M(X+Y-Z)(i+1,i-1)=-1$ holds, then it  increases the
switch-distance of the current $\widehat M$ by at most one (since it can be
directly back-swapped). The result of the second swap (after which the previous
problem is just solved automatically), together with our further strategy is
shown below:\\

\noindent
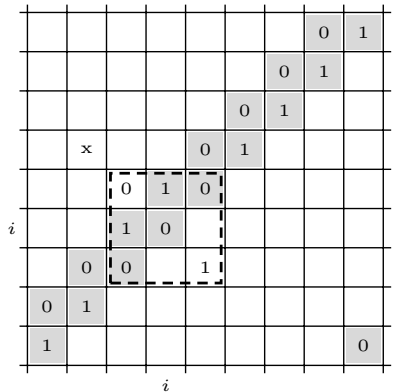
\begin{minipage}{6.5cm}
\scalebox{1.3}{\proba{
\begin{tikzpicture}[scale=0.4]
     \hh 40{$\scriptsize i$} \hh 04{ $\scriptsize i$} \kalap {9}
     \hh531    
    \hh350    \rhh330    \rhh550    \rhh341    \rhh451    \rhh440
    \draw[thick,\haha]  (3.1,3.1) rectangle (5.9,5.9);
    \hh26x\rhh{9}{1}{0}
\end{tikzpicture}
	}}
\end{minipage}
\quad
\begin{minipage}{5.2cm}
Here we distinguish between two cases, according to the value
$$F_Z(i-2,i+2)=x.$$  
This value can be $x=3$ or $x=0.$ 
\end{minipage}
In the case of $x=3$ 
we perform one more swap, which results in a subproblem with a friendly path (the swap shown on the left side of Figure~\ref{fig4}, while the right hand side indicates the two new subproblems):

\noindent
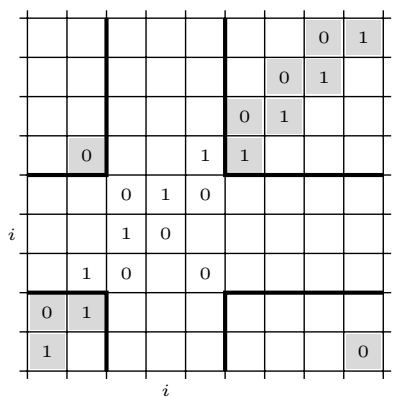
\begin{figure}[h]
\caption{The case of $x=3$}\label{fig4}
\scalebox{1.3}{\proba{
\begin{tikzpicture}[scale=0.4]
     \hh 40{$\scriptsize i$} \hh 04{ $\scriptsize i$} \kalap {9}
      \hh531
     \hh350    \rhh330    \rhh550    \rhh341    \rhh451  \rhh440
     \hh 263 
    \draw[thick,\haha]  (2.1,3.1) rectangle (5.9,6.9); \rhh{9}{1}{0}
\end{tikzpicture}
\ \raisebox{2cm}{$\Rightarrow$}

\begin{tikzpicture}[scale=0.4]
     \hh 40{$\scriptsize i$} \hh 04{ $\scriptsize i$} \kalap {9}
    \hh530
    \hh350    \rhh330    \rhh550    \rhh341    \rhh451   \rhh440    \rhh231    \rhh561
    \hh 262 
    \draw[thick,\haha]  (3.1,3.1) rectangle (5.9,5.9);
    \draw[very thick]  (1,3)--(3,3)--(3,1);
    \draw[very thick]  (1,6)--(3,6)--(3,10);
    \draw[very thick]  (10,6)--(6,6)--(6,10);
    \draw[very thick]  (6,1)--(6,3)--(10,3);\rhh{9}{1}{0}
\end{tikzpicture}
	}}
\end{figure}
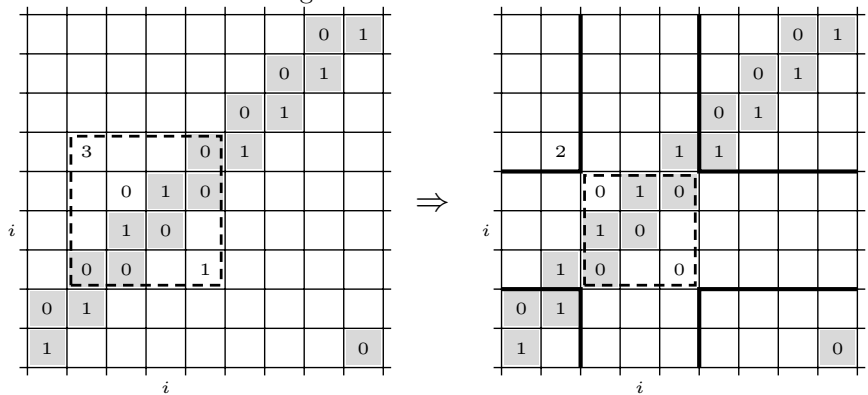

\noindent
\begin{minipage}{6cm}
\noindent \scalebox{1.3}{
\begin{tikzpicture}[scale=0.4]
     \hh 40{$\scriptsize i$} \hh 04{ $\scriptsize i$} \kalap {9}
    \hh530
    \hh350    \hh330    \hh550    \hh341    \hh451   \hh440    \hh231    \hh561
    \rhh 260 
    \draw[very thick]  (1,3)--(3,3)--(3,1);
    \draw[very thick]  (1,6)--(3,6)--(3,10);
    \draw[very thick]  (10,6)--(6,6)--(6,10);
    \draw[very thick]  (6,1)--(6,3)--(10,3);\rhh{9}{1}{0}
\end{tikzpicture}
}
\end{minipage}
\quad
\begin{minipage}{5.5cm}
Here we redrawn the RHS of Figure~\ref{fig4} to show that one subproblem is already solved: the cross shaped midsection  is in the required state: the main-diagonal (within the midsection) contains all $0$'s, while each entry in the small-diagonal is $1$ and the off-diagonal positions are in their original states. Here is nothing more to do.
\end{minipage}

\noindent  The second subproblem (indicated with the thick black lines, the four pieces fit together to a square matrix) is in the right form for further processing. The position $(2,6)$ changed into $0$ since we described the subproblem in the language of  (the now smaller) $F_Z$: the positions in the small-diagonal depend on the edges of $Z$ only.

We will process this second  subproblem along the up-line, containing position $(i-2,i+2)$, so the only currently improper entry will have the right value at the end of the next recursion step (that is it will be swapped back to its original value). (Here we can argue the same way as in Case 1.)

As it happened before  $\widehat M(X+Y-Z)$ may contain $2$ at the position
$(i-2,i+2)$. Again this  increases the switch-distance  by at most one, since
the positions $(i-2,i-1), (i+1,i-1)$ and $(i+1,i+2)$ are not in our subproblem.

\bigskip\noindent Finally it can happen, that $x=0.$ Then we can define the
following subproblem:

\noindent
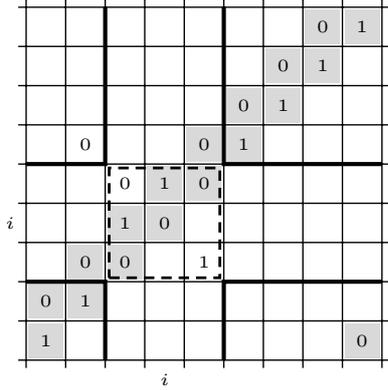
\begin{minipage}{6.5cm}
\noindent
\scalebox{1.3}{
\begin{tikzpicture}[scale=0.4]
    \hh 40{$\scriptsize i$} \hh 04{ $\scriptsize i$} \kalap {9}
    \hh531    \hh350    \rhh330    \rhh550    \rhh341    \rhh451
    \rhh440    \hh260
    \draw[thick,\haha]  (3.1,3.1) rectangle (5.9,5.9);
    \draw[very thick]  (1,3)--(3,3)--(3,1);
    \draw[very thick]  (1,6)--(3,6)--(3,10);
    \draw[very thick]  (10,6)--(6,6)--(6,10);
    \draw[very thick]  (6,1)--(6,3)--(10,3);\rhh{9}{1}{0}
\end{tikzpicture}
	}
\end{minipage}
\begin{minipage}{5.2cm}
\noindent This figure shows the new subproblem (indicating with thick black
lines) is in the right form (for further processing)  again. We will process the subproblem along the up-line, containing position $(i-2,i+2)$ (so the only currently improper entry
will have the right value at the end of the next recursion step).
\end{minipage}

\noindent
\begin{minipage}{8cm}
Here, again, we may confront the fact, that $\widehat M(X+Y-Z)(i+1,i-1)=-1.$ Then we
should consider the alternating cycle shown in the figure. All elements of the
cycle, except $(i+1,i-1)$, are in the main- and small-diagonal, therefore along
this cycle we can swap that entry into range within a small number (say
$\delta$) of steps. This will increase the switch-distance of $\widehat M$ by at most $\delta.$
\end{minipage}
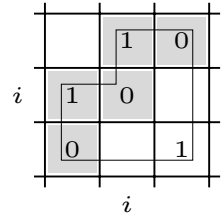
\begin{minipage}{3cm}
\scalebox{1.8}{\proba{
\begin{tikzpicture}[scale=0.4]
    \kalap {3}
    \hh 02{$\scriptsize i$} \hh 20{$\scriptsize  i$}
    \rhh110 \rhh220 \rhh330 \rhh121 \rhh231 \hh311
    \draw [very thin] (1.3,1.3) -- (1.3,2.7) -- (2.3,2.7) -- (2.3,3.7) --
(3.7,3.7) -- (3.7,1.3) -- (1.3,1.3);
\end{tikzpicture}
	}}
\end{minipage}

\smallskip
\noindent We run the  first recursion on the subproblem along the $i$th up-line, therefore the sub-subproblem with friendly path will contain the position $(i-2, i+2)$. Therefore  when we finish the first recursion, our matrix $F_G$ will be in the following form: (the figure on the left):

\noindent\scalebox{1.3}{\proba
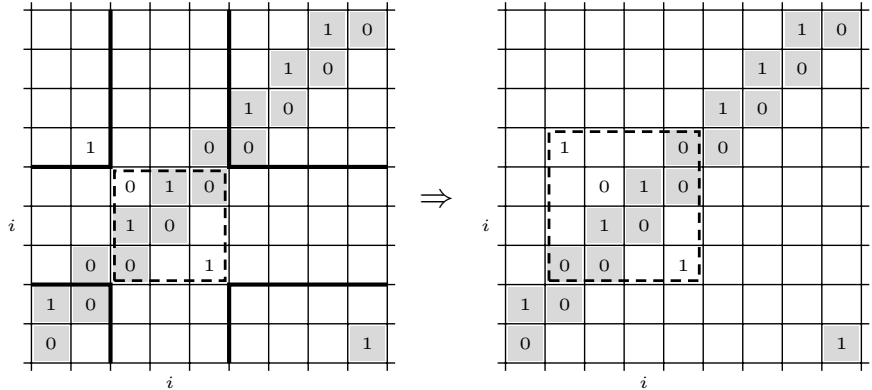
{
\begin{tikzpicture}[scale=0.4]
     \hh 40{$\scriptsize i$} \hh 04{$ \scriptsize i$} \kalap {9}
    \hh531 \hh350 \rhh330 \rhh550 \rhh341 \rhh451 \rhh440 \hh261
    \rhh110 \rhh220 \rhh121 \rhh660 \rhh671
    \rhh770  \rhh 880 \rhh 990 \rhh 781  \rhh 891
    \draw[thick,\haha]  (3.1,3.1) rectangle (5.9,5.9);
    \draw[very thick]  (1,3)--(3,3)--(3,1);
    \draw[very thick]  (1,6)--(3,6)--(3,10);
    \draw[very thick]  (10,6)--(6,6)--(6,10);
    \draw[very thick]  (6,1)--(6,3)--(10,3);\rhh{9}{1}{1}
\end{tikzpicture}
\ \raisebox{2cm}{$\Rightarrow$}
\begin{tikzpicture}[scale=0.4]
     \hh 40{$\scriptsize i$} \hh 04{$ \scriptsize i$} \kalap {9}
    \hh531 \hh350 \rhh220 \rhh121 \rhh660 \rhh671
    \rhh330 \rhh550 \rhh341 \rhh451 \rhh440 \hh261 \rhh110
    \rhh770 \rhh 880 \rhh 990 \rhh 781  \rhh 891
    \draw[thick,\haha]  (2.1,3.1) rectangle (5.9,6.9);\rhh{9}{1}{1}
\end{tikzpicture}
	}}

\medskip
\noindent We have seen how one can handle the switch-distance of our matrix, if
position $(i+1,i-1)$ is problematic with $\widehat M(i+1,i-1)=-1$ but position $(i-2,i+2)$ is correct. On the other hand if $\widehat M(X+Y-Z)(i-2,i+2)=-1$  then the swap on the positions $(i-2,i-1),(i+1,i+2); (i+1,i-1),(i-2,i+2)$ changes both $(i+1,i-1)$ and $(i-2,i+2)$ into $0.$ For $(i+1,i-1)$ that was the original type - so it cannot be wrong in $\widehat M.$

\noindent
\begin{minipage}{8cm}
\noindent After that we perform the swap on the positions $(i-2,i-1),(i+1,i+2);
(i+1, i-1),(i-2,i+2)$ (these are the corners of the dashed square in the figure
on the upper right). The result is shown to the right:
\end{minipage}
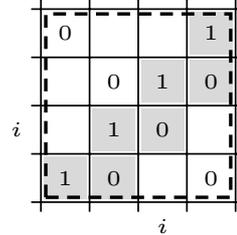
\begin{minipage}{3cm}
\scalebox{1.6}{\proba{
\begin{tikzpicture}[scale=0.4]
     \hh 30{$\scriptsize i$}  \hh 02{$\scriptsize i$} \kalap {4}
    \hh 120 \hh 12{} \hh 230 \hh 340 \hh 34{} \rhh210 \rhh320 \rhh430
    \hh140
    \hh410
    \draw [thick,\haha] (1.1,1.1) rectangle (4.9,4.9);
\end{tikzpicture}
	}}
\end{minipage}

\noindent This completes our handling on the First Possibility, that is when
for our $i$ we have the  value $j'=1.$ Now we turn to the other (and probably
more common) configuration:

\medskip\noindent{\bf Second Possibility:} We have $j'\ge 2$. Unfortunately, the
situation can be more complicated in this case due to the possible switch-distances of
$\widehat M.$ We overcome this problem by showing at first the general structure of
the process, and later we give the necessary fine-tuning to ensure the bounded
switch-distance. (Recall again, that the bounded switch distance is necessary to have a good upper bound on the number of different matrices $\widehat M$ appearing along the algorithm.)

\bigskip\noindent In our current alternating cycle (lying in the symmetric difference of $G$ and $G'$) there is no friendly path, therefore there is a Steinhaus set $T$ in  $F_G$. Now fix a particular $i$ and assume that the $j'$ corresponding to this $i$ is $\ge 2.$ We should distinguish between two cases: where the down-line starts with the value $t=3$
or with $t=0.$  

\smallskip\noindent {\bf Case 1: $t=3$} The first figure below shows the structure of matrix $F_G.$ The dashed square is the first subproblem to deal with, while the thick black lines indicate the second subproblem. However, before we start the processing the subproblems, we have to perform a swap. The corners of the thin black square shows the positions of the swap.

\medskip
\noindent\scalebox{1.3}{\proba
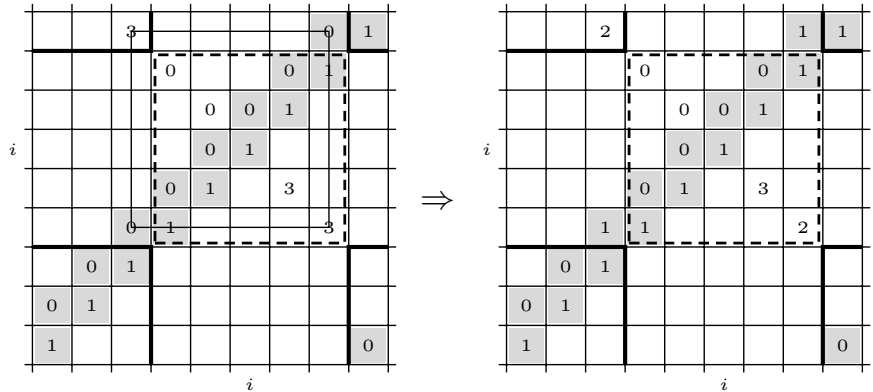
{
\begin{tikzpicture}[scale=0.4]
    \hh 60{$\scriptsize i$}  \hh 06{$\scriptsize i$} \kalap {9}
    \hh 753  
    \hh 843 \hh570 \hh 480 \hh 393  
    \draw (3.5,4.5) -- (3.5,9.5) -- (8.5,9.5) -- (8.5,4.5) -- (3.5,4.5);
    \draw[thick,\haha]  (4.1,4.1) rectangle (8.9,8.9);
    \draw[very thick]  (1,4)--(4,4)--(4,1);
    \draw[very thick]  (1,9)--(4,9)--(4,10);
    \draw[very thick]  (10,9)--(9,9)--(9,10);
    \draw[very thick]  (9,1)--(9,4)--(10,4); \rhh{9}{1}{0}
\end{tikzpicture}
\ \raisebox{2cm}{$\Rightarrow$}
\begin{tikzpicture}[scale=0.4]
    \hh 60{$\scriptsize i$}  \hh 06{$\scriptsize i$} \kalap {9}
    \hh 753 \hh 842 \hh570 \hh 480 \hh392 \rhh 341 \rhh891
    \draw[thick,\haha]  (4.1,4.1) rectangle (8.9,8.9);
    \draw[very thick]  (1,4)--(4,4)--(4,1);
    \draw[very thick]  (1,9)--(4,9)--(4,10);
    \draw[very thick]  (10,9)--(9,9)--(9,10);
    \draw[very thick]  (9,1)--(9,4)--(10,4);\rhh{9}{1}{0}
\end{tikzpicture}	
	}}

\noindent After this swap (see the figure above, right.), the first subproblem (indicated by the dashed square) is in the right form. Indeed, the left figure below shows  the two separate subproblems.

\noindent \scalebox{1.3}
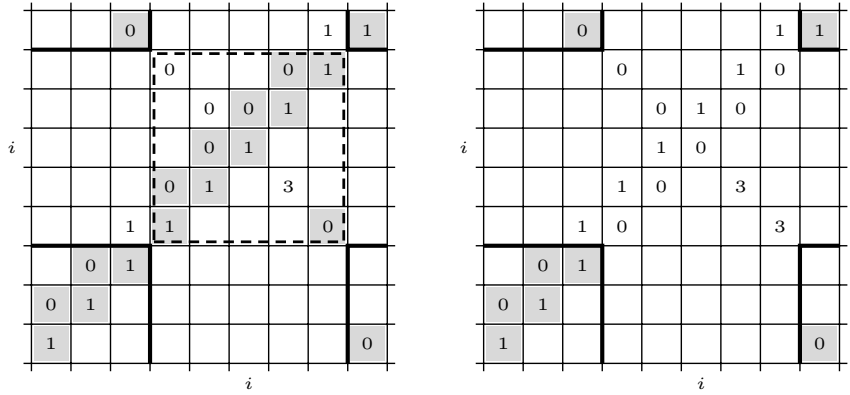
{
\begin{tikzpicture}[scale=0.4]
    \hh 60{$\scriptsize i$}  \hh 06{$\scriptsize i$} \kalap {9}
    \hh 753 \rhh 840 \hh570 \hh 480 \rhh390 \rhh 341 \rhh891
    \hh 341 \hh891
    \draw[thick,\haha]  (4.1,4.1) rectangle (8.9,8.9);
    \draw[very thick]  (1,4)--(4,4)--(4,1);
    \draw[very thick]  (1,9)--(4,9)--(4,10);
    \draw[very thick]  (10,9)--(9,9)--(9,10);
    \draw[very thick]  (9,1)--(9,4)--(10,4);\rhh{9}{1}{0}
\end{tikzpicture}	
     } \quad
\noindent\scalebox{1.3}{
\begin{tikzpicture}[scale=0.4]
    \hh 60{$\scriptsize i$}  \hh 06{$\scriptsize i$} \kalap {9}
    \hh 753 \hh 843 \hh570 \hh 480 \rhh390 \hh 341 \hh891
    \hh 440 \hh 550 \hh 660 \hh 770 \hh 880
    \hh 451 \hh 561 \hh 671 \hh 781
    \draw[very thick]  (1,4)--(4,4)--(4,1);
    \draw[very thick]  (1,9)--(4,9)--(4,10);
    \draw[very thick]  (10,9)--(9,9)--(9,10);
    \draw[very thick]  (9,1)--(9,4)--(10,4);\rhh{9}{1}{0}
\end{tikzpicture}	
	}

\noindent Finishing the first subproblem, we have the  $F_Z$ matrix (above, right). As it can be seen, after the first phase, all entries in the midsection are in their required types: the small-diagonal consists of $1$s (including position $(i+2, i-2)$ which in that way is back to its original type), while the main-diagonal consists of only $0$s.

The second subproblem (indicated by the thick black lines) in the right form now to process
(including position $(i-3,i+3)$ which is sitting on the small-diagonal).

\noindent After completing the solution of the black subproblem, all entries in
the matrix will be in exactly the required type. We start processing the black
subproblem  on the $i$th up-line, therefore the actual types of positions
$(i-3,i+3)$ and $(i+2, i-2)$ can be described as follows: Position $(i+2,i-2)$
has opposite type after the very first  swap, then while processing the dashed
subproblem it may change between $0$ and $1$. Finishing the dashed subproblem,
it will be in the same type as it starts.

Position $(i-3,i+3)$ will be in type $0$ all the way in the dashed phase, while
within the black phase it will change between $1$ and $0.$ At the end, as we
already mentioned, is $1.$

\bigskip \noindent {\bf Case 2: $t=0$} The first figure below shows the
structure of matrix $F_G.$ The dashed square is the first subproblem to deal
with, while the thick black lines indicate the second subproblem. They can
process without any preprocessing.
\bigskip
\noindent\scalebox{1.3}{\proba
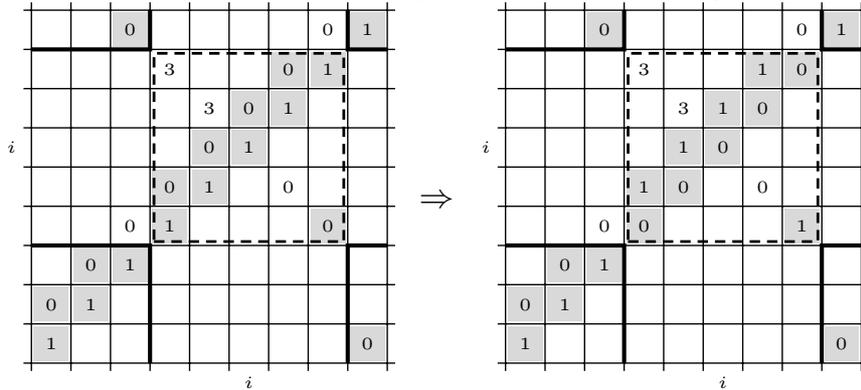
{
\begin{tikzpicture}[scale=0.4]
    \hh 60{$\scriptsize i$}  \hh 06{$\scriptsize i$} \kalap {9}
    \hh 750 \rhh 840 \hh573 \hh 483 \rhh390 \hh340 \hh890
    \draw[thick,\haha]  (4.1,4.1) rectangle (8.9,8.9);
    \draw[very thick]  (1,4)--(4,4)--(4,1);
    \draw[very thick]  (1,9)--(4,9)--(4,10);
    \draw[very thick]  (10,9)--(9,9)--(9,10);
    \draw[very thick]  (9,1)--(9,4)--(10,4);\rhh{9}{1}{0}
\end{tikzpicture}
\ \raisebox{2cm}{$\Rightarrow$}
\begin{tikzpicture}[scale=0.4]
    \hh 60{$\scriptsize i$}  \hh 06{$\scriptsize i$} \kalap {9}
    \hh 750 \rhh 841 \hh573 \hh 483 \rhh390
    \rhh 440 \rhh 550 \rhh 660 \rhh 770 \rhh 880  \hh340 \hh890
    \rhh 451 \rhh 561 \rhh 671 \rhh 781
    \draw[thick,\haha]  (4.1,4.1) rectangle (8.9,8.9);
    \draw[very thick]  (1,4)--(4,4)--(4,1);
    \draw[very thick]  (1,9)--(4,9)--(4,10);
    \draw[very thick]  (10,9)--(9,9)--(9,10);
    \draw[very thick]  (9,1)--(9,4)--(10,4);\rhh{9}{1}{0}
\end{tikzpicture}	
	}}

\noindent At the end everything will be in the right type, except the four positions, showed by the thin black square (below, left side). We can finish the process with that  swap.

\smallskip
\noindent\scalebox{1.3}
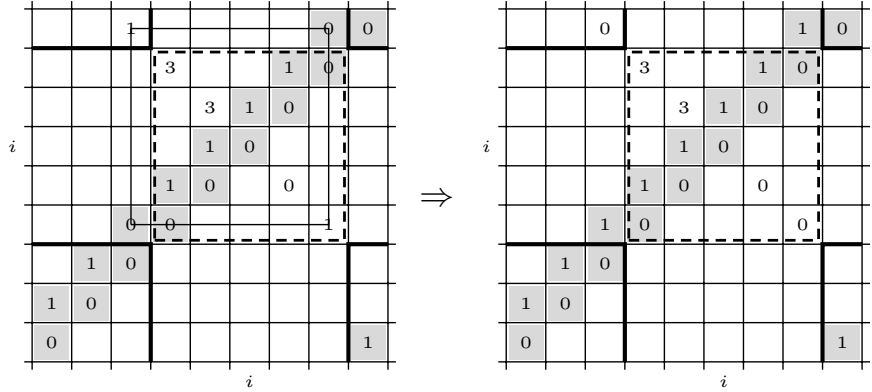
{
\begin{tikzpicture}[scale=0.4]
    \hh 60{$\scriptsize i$}  \hh 06{$\scriptsize i$} \kalap {9}
    \hh 750 \hh 841 \hh573 \hh 483 \hh391
    \rhh 440 \rhh 550 \rhh 660 \rhh 770 \rhh 880
    \rhh 451 \rhh 561 \rhh 671 \rhh 781
    \rhh 110 \rhh220 \rhh 330 \rhh990
    \rhh 121 \rhh 231 \rhh911
    \draw (3.5,4.5) -- (3.5,9.5) -- (8.5,9.5) -- (8.5,4.5) -- (3.5,4.5);
    \draw [thick,\haha] (4.1,4.1) rectangle (8.9,8.9);
    \draw[very thick]  (1,4)--(4,4)--(4,1);
    \draw[very thick]  (1,9)--(4,9)--(4,10);
    \draw[very thick]  (10,9)--(9,9)--(9,10);
    \draw[very thick]  (9,1)--(9,4)--(10,4);
\end{tikzpicture}	
\ \raisebox{2cm}{$\Rightarrow$}
\begin{tikzpicture}[scale=0.4]
    \hh 60{$\scriptsize i$}  \hh 06{$\scriptsize i$} \kalap {9}
    \hh 750 \hh 840 \hh573 \hh 483 \hh390
    \rhh 440 \rhh 550 \rhh 660 \rhh 770 \rhh 880
    \rhh 451 \rhh 561 \rhh 671 \rhh 781
    \rhh 110 \rhh220 \rhh 330 \rhh990
    \rhh 121 \rhh 231 \rhh911
    \rhh 341 \rhh 891
    \draw [thick,\haha] (4.1,4.1) rectangle (8.9,8.9);
    \draw[very thick]  (1,4)--(4,4)--(4,1);
    \draw[very thick]  (1,9)--(4,9)--(4,10);
    \draw[very thick]  (10,9)--(9,9)--(9,10);
    \draw[very thick]  (9,1)--(9,4)--(10,4);
\end{tikzpicture}	
	}

\bigskip\noindent While the overall structure of our plan is clear, we may meet
problems along this procedure. The reason is that we must be able to control the
switch-distance of our $\widehat M(X+Y-Z)$ (we will use here simply $\widehat M$) from
the adjacency matrix of some realization. There  are two neuralgic points: both
the positions $(i+j,i-j)$ and $(i-j',i+j')$ may contain \old{the} $-1$, or both
may contain $2$.  When we start a new subproblem, then their types always
provide a suitable switch for the control (as it was seen before). However, when
we proceed along our subproblem, then it can happen that one of the problematic
positions changes its value, while the other does not. But in this case the
switch which was previously available is not useable anymore. Next we
describe how we can fine tune our procedure to avoid this trap.

As we know the first subproblem contains a friendly path (by Corollary \ref{lm:van-friendly}), and for easier reference let call its problematic position $P_1.$ We also know that second
subproblem contains a problematic position, $P_2,$ and probably we have to divide this subproblem into two smaller ones. If so, then the first of them becomes the new second subproblem, which contains $P_2$ and possesses a friendly path, while the third  subproblem contains another problematic position, $P_3.$

\medskip\goodbreak
\noindent{\bf Fine tuning:}
\begin{enumerate}
 \item We begin our swap sequence along the first subproblem but we stop just
before we face the swap which changes the value of $P_1.$
 \item Next we continue with the swap sequence of the second problem and we stop
before we should perform a swap on $P_2.$
 \item Now we finish the swap sequence of the first subproblem.
 \item After that we focus on the second subproblem. Dealing  effectively with
this, we need to prepare the third subproblem similarly as we did with the
second one, when we were working on the first one. Therefore we begin the swap
sequence of the third subproblem but we stop it before the first swap would be
carried out on $P_3.$
 \item And if now we just rename our two active subproblems as first and second
subproblem, we are back to a situation, which is equivalent to the beginning of
the third stage.
\end{enumerate}
\noindent
Along this algorithm, at each point we have two "active" subproblems. When a subproblem has a friendly path, then along this path we define the necessary swap sequence (as described in Subsection \ref{sec:friendly}) and we have an upper bound on its length). When the subproblem is without a friendly path, then we divide it into two, and one (or both) of them have a friendly path, etc. The sum of the sizes of the subproblems is at most the size of the original cycle. Finally we put together the final swap sequence from these swap sequences and  some short sequences we get from the (sometimes) necessary preprocessing. Finally, since we have bounded switch distances all along (one or two at preprocessing stages, and those given in Subsection  \ref{sec:friendly}), therefore all together we have a good control of the overall number of used $\widehat M$'s. \qed

\section{Acknowledgement}

The authors would like to thank to the anonymous referee, whose comments and suggestions improved the manuscript significantly. We are most grateful to Catherine Greenhill for her tremendous help to prepare this manuscript.

\end{document}